\documentclass[12pt]{elsarticle}

\usepackage{latexsym,amsfonts,amsmath,theorem,amssymb}
\usepackage{graphicx}

\def\be{\begin{equation}}
\def\bea{\begin{eqnarray*}}
\def\ee{\end{equation}}
\def\eea{\end{eqnarray*}}
\def\ba{\begin{array}}
\def\ea{\end{array}}
\def\bi{\begin{itemize}}
\def\ei{\end{itemize}}
\newtheorem{theo}{Theorem}
\newtheorem{lem}{Lemma}

\def\tP{\tilde{P}}

\def\tH{\tilde{H}}

\def\th{\tilde{h}}

\def\mF{\mathcal{F}}

\journal{XXX}

\begin{document}
\begin{frontmatter}

\title{Multivariate Haar systems in Besov function spaces}
\author[bonn]{Peter Oswald}
\ead{agp.oswald@gmail.com}

%\thanks{}
\address[bonn]{Institute for Numerical Simulation (INS),
University of Bonn,
Wegelerstr. 6-8,
D-53211 Bonn
}
%\thanks{}
\date{}
%\maketitle
\begin{abstract}
We determine all cases for which the $d$-dimensional Haar wavelet system $H^d$ on the unit cube $I^d$ is a conditional or unconditional Schauder basis in the classical isotropic Besov function spaces ${B}_{p,q,1}^s(I^d)$, $0<p,q<\infty$, $0\le s < 1/p$, defined in terms of first-order $L_p$ moduli of smoothness. We obtain similar results for the tensor-product Haar system $\tH^d$, and characterize the parameter range for which the dual of ${B}_{p,q,1}^s(I^d)$ is trivial for $0<p<1$.
\end{abstract}
\end{frontmatter}

\section{Introduction}\label{sec1} The univariate Haar system $H:=\{h_m\}_{m\in\mathbb{N}}$ was one of the first examples 
of a Schauder basis in some classical function spaces on the unit interval $I:=[0,1]$, see \cite{Ci1977}, \cite[Section III]{KaSa1989}, and \cite[Section 2.1]{Tr2010} for a review of the early history of the Haar system as basis in function spaces. Meantime the existence of Schauder bases in function spaces of Besov-Hardy-Sobolev type has been established in most cases, see \cite{Tr2006} for a recent exposition. Early on, a major step was taken by Ciesielski  and co-workers \cite{CiDo1972,Ci1975,Ci1977,CiFi1982,CiFi1983a,CiFi1983b} who constructed families of spline systems generalizing the classical Haar, Faber, and Franklin systems, and established their basis properties in Lebesgue-Sobolev spaces over $d$-dimensional cubes and smooth manifolds for $1\le p \le \infty$. For distributional Besov spaces $B_{p,q}^s(\mathbb{R}^d)$ and Triebel-Lizorkin spaces $F_{p,q}^s(\mathbb{R}^d)$ with $0<p,q<\infty$, $s\in \mathbb{R}$, wavelet systems provide examples of unconditional Schauder bases. Results in these directions and their generalization to spaces on domains in $\mathbb{R}^d$ are presented in \cite{Tr2008,Tr2010}. Needless to say that not all quasi-Banach function spaces possess nice basis properties. E.g., $L_1(I)$ does not possess an unconditional Schauder basis, see \cite[Theorem II.13]{KaSa1989}, while the quasi-Banach space $L_p(I)$, $0<p<1$, cannot have Schauder bases at all since its dual $L_p(I)'=\{0\}$ is trivial.

%The question if a concrete system is a basis in a concrete function space or scale of function spaces is also of interest. In particular, low-order spline systems such as Faber- and Franklin-type systems on polyhedral domains in $\mathbb{R}^d$ often appear in investigations on numerical integration and multiscale finite element methods, see \cite{Os1994,Tr2010}.   

In this paper, we deal with the multivariate anisotropic tensor-product Haar system
\be\label{HTPS}
\tH^d = \underbrace{H\otimes \ldots \otimes H}_{d \mbox{ times}},\qquad d\ge 1,
\ee
and its isotropic counterpart $H^d$ on the unit cube $I^d\subset \mathbb{R}^d$ (the latter is called \emph{Haar wavelet system} in \cite{Tr2010}), and consider their Schauder basis properties in the Besov spaces ${B}_{p,q,1}^s(I^d)\subset L_p(I^d)$. These function spaces are classically defined in terms of first-order $L_p$ moduli of smoothness (detailed definitions are given in the next section), and coincide with their distributional counterparts
${B}_{p,q}^s(I^d)$ only under some restrictions on $p,q,s$. We mostly concentrate on the parameter range
\be\label{PR}
0<p,q<\infty,\qquad 0 < s < 1/p.
\ee
With the exception of the special case $s=0$, this is the maximal range of parameters for which ${B}_{p,q,1}^s(I^d)$ is a separable quasi-Banach space \emph{and} contains the Haar systems $\tH^d$ and $H^d$. Moreover, for this parameter range ${B}_{p,q,1}^s(I^d)$ admits a characterization in terms of best $L_p$-approximations with piecewise constant functions on dyadic partitions of $I^d$ which is used in the proofs. Our main result for the Haar wavelet system $H^d$ is the following theorem.
\begin{theo}\label{theo1} Assume (\ref{PR}).\\
a) If $1\le p<\infty$ then the Haar wavelet system $H^d$ is an unconditional Schauder basis for ${B}_{p,q,1}^s(I^d)$ for all parameters
$0< s <1/p$, $0<q<\infty$ of interest. \\
b) Let  $0<p<1$. Then $H^d$ is an unconditional Schauder basis for ${B}_{p,q,1}^s(I^d)$ if and only if $d(1/p-1)<s<1/p$, $0<q<\infty$.
If $s=d(1/p-1)$, $0<q\le p$ then $H^d$ is a Schauder basis for ${B}_{p,q,1}^s(I^d)$ but not unconditional. In all other cases, $H^d$ does not possess the Schauder basis property for ${B}_{p,q,1}^s(I^d)$.
\end{theo}

The statements about unconditionality of $H^d$ in Theorem \ref{theo1} are proved by giving a characterization of ${B}_{p,q,1}^s(I^d)$ in terms of Haar coefficients. We refer to Theorem \ref{theo3} in Section \ref{sec32} which states conditions under which  ${B}_{p,q,1}^s(I^d)$ is isomorphic to some weighted $\ell_q(\ell_p)$ sequence space. 
The exceptional case $s=0$ is covered by Theorem \ref{theo5} in Section \ref{sec53}.
We also have
\begin{theo}\label{theo2} If $0<p<1$ in (\ref{PR}) then the Besov space ${B}_{p,q,1}^s(I^d)$ does not possess nontrivial bounded linear functionals, i.e.,
	${B}_{p,q,1}^s(I^d)'=\{0\}$, if and only if $s<d(1/p-1)$, $0<q<\infty$ or $s=d(1/p-1)$, $1<q<\infty$. 
\end{theo}

Consequently, for these parameters ${B}_{p,q,1}^s(I^d)$ does not have a Schauder basis at all which is a  stronger statement than proving that the Haar wavelet system $H^d$ fails to be a Schauder basis. For the parameters
\be\label{PR1}
0<p <q \le 1,\qquad s = d(1/p-1),
\ee
where according to Theorem \ref{theo1} $H^d$ is not a Schauder basis in ${B}_{p,q,1}^s(I^d)$, we have the continuous embedding ${B}_{p,q,1}^s(I^d)\subset L_1(I^d)$, and thus $L_\infty(I^d)\subset {B}_{p,q,1}^s(I^d)'$. We do not know if the spaces ${B}_{p,q,1}^s(I^d)$ satisfying (\ref{PR1}) possess (unconditional) Schauder bases at all.

As Schauder basis in Besov spaces ${B}_{p,q,1}^s(I^d)$, the tensor-product Haar system $\tH^d$
behaves the same as $H^d$ for $1\le p < \infty$ but fails completely for the parameter arnge $p<1$. This is the essence of Theorem \ref{theo4} in Section  \ref{sec52}.

Let us comment on known results that motivated this study. For $d>1$, the two Haar systems $\tH^d$ and
$H^d$ have formally been introduced in the 1970-ies, see \cite{Ci1977,CiFi1983b}, the system $H^d$ implicitly appeared already in \cite{Tr1973}. In the univariate case $d=1$, part a) of Theorem \ref{theo1} has essentially been established by Triebel \cite{Tr1973} and Ropela \cite{Ro1976} under the restriction $q\ge 1$. Extensions to $d\ge 1$ are due to Ciesielski \cite{Ci1977} and Triebel. The results of Triebel who worked in the framework of distributional Besov spaces $B_{p,q}^s(\mathbb{R}^d)$ and $B_{p,q}^s(I^d)$ are summarized in \cite[Section 2]{Tr2010}. Starting from \cite{Tr1978}, Triebel considered the parameter values $0<p,q<\infty$, $-\infty<s<\infty$, and proved the unconditionality of $H^d$ in 
$B_{p,q}^s(I^d)$ for the parameter range 
$$%\be\label{PR21}
\max(d(1/p-1),1/p-1) < s < \min(1,1/p), \quad 0 < p,q <\infty.
$$%\ee
This is the essence of Theorem 2.13 (i) ($d=1$) and Theorem 2.26 (i) ($d>1$) in \cite{Tr2010}.
Note that for the range $0 <p,q\le \infty$ we have 
\be\label{PR31}
B_{p,q}^s(I^d) = B_{p,q,1}^s(I^d) \quad \Longleftrightarrow \quad \max(0,d(1/p-1)) < s < \min(1,1/p),
\ee
i.e., for the parameters in (\ref{PR31}), the scales $B_{p,q}^s(I^d)$ and $B_{p,q,1}^s(I^d)$ coincide up to equivalent quasi-norms. Consequently, with the exception of the range $1\le s <1/p$ for $0<p<1$,  the unconditionality of $H^d$ in the spaces $B_{p,q,1}^s(I^d)$ essentially follows from Triebel's results for $B_{p,q}^s(I^d)$ for all cases stated in Theorem \ref{theo1} and Theorem \ref{theo3}. However, we give a direct proof for $B_{p,q,1}^s(I^d)$
using its characterization by piecewise constant best $L_p$ approximations on dyadic partitions. 

Triebel \cite{Tr1978} also established that outside the \emph{closure} of the parameter range (\ref{PR21}), the Haar wavelet system $H^d$ is not a Schauder basis in the spaces $B_{p,q}^s(I^d)$. The boundary cases remained unsettled until recently when Garrig\'os, Seeger, and Ullrich  dealt in a series of papers \cite{GSU2018,GSU2019a,GSU2019b,SeUl2017a,SeUl2017b} with the open cases for both the
distributional $B_{p,q}^s$ and $F_{p,q}^s$ scales. In particular,
\cite{GSU2019a} provides complete answers concerning the Schauder basis properties of the Haar wavelet systems in
$B_{p,q}^s(\mathbb{R}^d)$ and $B_{p,q}^s(I^d)$. They also established a subtle difference between the cases $\mathbb{R}^d$ and $I^d$ for the critical smoothness  parameter $s=d(1/p-1)$, $0<p<1$, and provided correct asymptotic estimates of the norms of partial sum projectors associated with $H^d$. We also mention the paper \cite{YSY2018} related to the questions considered in this paper, where the authors study necessary and sufficient conditions on the parameters $p,q,s,\tau$ for which the map
$f \to (f,\chi_{I^d})_{L_2}=\int_{I^d} f\,dx$ extends to a bounded linear functional on Besov-Morrey-Campanato-type spaces $B_{p,q}^{s,\tau}(\mathbb{R}^d)$.

Independently, for $0<p<1$ the Schauder basis property of $H^d$ in $B_{p,q}^s(I^d)$ and $B_{p,q,1}^s(I^d)$
spaces was also considered by the author in \cite{Os2018}. This paper expanded on \cite{Os1981}, where a partial result was stated for $d=1$, namely that the univariate Haar system $H$ forms a Schauder basis in $B_{p,q,1}^{s}(I)$  for the critical smoothness parameter $s=1/p-1$ if $0<q\le p<1$ (see the remark at the end of \cite{Os1981}). In \cite{Os1981}, it was also established that $B_{p,q,1}^{s}(I)$ has a trivial dual if $0<s<1/p-1$, $0<p<1$, $0<q<\infty$.

The paper is organized as follows. In Section \ref{sec2}, we give the necessary definitions and state auxiliary results on Besov spaces and on piecewise constant $L_p$-approximation with respect to dyadic partitions. Section \ref{sec3} deals with the proof of the sufficiency of the conditions on the parameters $p,q,s$ appearing in the main results formulated in Theorems \ref{theo1} and \ref{theo3}. The necessity of these conditions and Theorem \ref{theo2} are dealt with in Section \ref{sec4}, where we construct specific counterexamples consisting of piecewise constant functions. Similar in spirit examples have already been used in \cite{Os1981}. We conclude in Section \ref{sec5} with some remarks on higher-order spline systems, analogous results for the tensor-product Haar system $\tH^d$, and the exceptional cases $s=0$ and $q=\infty$.

\section{Definitions and auxiliary results}\label{sec2}

\subsection{Haar systems}\label{sec21}
Recall first the definition of the univariate $L_\infty$ normalized Haar functions. By $\chi_\Omega$ we denote the characteristic function of a Lebesgue measurable set $\Omega\subset \mathbb{R}^d$, and
by $\Delta_{k,i}:=[(i-1)2^{-k},i2^{-k})$ the univariate dyadic intervals of length $2^{-k}$, $k\in \mathbb{Z}_+$, $i\in\mathbb{Z}$. Then the univariate Haar system
$H=\{h_m\}_{m\in\mathbb{N}}$ on $I:=[0,1]$ is given by
$h_1=\chi_I$, and
$$
h_{2^{k-1}+i}=\chi_{\Delta_{k,2i-1}}-\chi_{\Delta_{k,2i}}, \qquad i=1,\ldots,2^{k-1}, \quad k\in\mathbb{N},
$$
Throughout the paper, we work with $L_\infty$-normalized Haar functions.
The Haar functions $h_m$ with $m\ge 2$ can also be indexed
by their supports, and identified with the shifts and dilates
of a single function, the Haar wavelet
$h_0:=\chi_{[0,1/2)}-\chi_{[1/2,1)}$.
Indeed,
$$
h_{\Delta_{k-1,i}}:=h_{2^{k-1}+i}=h_0(2^{k-1} \cdot-i),\qquad i=1,\ldots,2^{k-1},\quad k\in\mathbb{N}.
$$
The above introduced enumeration of the Haar functions $h_m$ is the natural ordering used in the literature,
however, one can also define $H$ as the union of dyadic blocks
$$
H=\cup_{k=0}^\infty H_k,\qquad H_0=\{h_1\},\quad H_k=\{h_{\Delta_{k-1,i}}:\,i=1,\ldots,2^{k-1}\},\quad k\in \mathbb{N},
$$
and allow for arbitrary orderings within each block $H_k$.  
Below, we will work with the multivariate counterparts of the spaces
$$
S_k=\mathrm{span}(\{h_m\}_{m=1}^{2^k})=\mathrm{span}(\{\chi_{\Delta_{k,i}}\}_{i=1}^{2^k}),\qquad k=0,1,\ldots,
$$
of piecewise constant functions with respect to the uniform dyadic partition $T_k=\{\Delta_{k,i}:\,i=1,\ldots,2^k\}$ of step-size $2^{-k}$ on the unit interval $I$. %, which we call for short \emph{dyadic step functions of level} $k$.

The Haar wavelet system 
\be\label{HWS}
H^d = \cup_{k=0}^\infty H_k^d
\ee
on the $d$-dimensional cube $I^d$, $d>1$, is defined in a blockwise fashion as follows. Let the partition $T_k^d$ be the set of all dyadic cubes of side-length $2^{-k}$ in $I^d$. Each cube in $T_k^d$ is the $d$-fold 
product of univariate $\Delta_{k,i}$, i.e., 
$$
T_k^d =\{ \Delta_{k,\mathbf{i}}:=\Delta_{k,i_1}\times \ldots\times \Delta_{k,i_d}:\;\mathbf{i}=(i_1,\ldots,i_d)\in \{1,\ldots,2^k\}^d\}.
$$
The set of all piecewise constant functions on $T_k^d$ is denoted by $S_k^d$.
With each $\Delta_{k-1,\mathbf{i}}\in T_{k-1}^d$, $\mathbf{i}\in \{1,\ldots,2^{k-1}\}^d$, $k\in \mathbb{N}$,
we associate the set $H^d_{k,\mathbf{i}}\subset S_k^d$ of $2^d-1$ multivariate Haar functions with support $\Delta_{k-1,\mathbf{i}}$, given by all possible tensor products
$$
\psi_{k,i_1}\otimes \psi_{k,i_2}\otimes \ldots \otimes\psi_{k,i_d},\qquad \psi_{k,i} = h_{\Delta_{k-1,i}}\mbox{ or } \chi_{\Delta_{k-1,i}}
$$
where at least one of the $\psi_{k,i_l}$ equals $h_{\Delta_{k-1,i_l}}$.
The blocks $H_k^d$ appearing in (\ref{HWS}) are given as follows:
The block $H_0^d$ is exceptional, and  consists of the single constant function $\chi_{I^d}$. The block $H_1^d$ coincides with 
$H^d_{1,\mathbf{1}}$ and consists of $2^d-1$ Haar functions, where $\mathbf{1}=(1,\ldots,1)$. For general $k\ge 2$, the block
$$
H_k^d := \cup_{\Delta_{k-1,\mathbf{i}}\in T_{k-1}^d} \, H^d_{k,\mathbf{i}}
$$
consists of $(2^d-1)2^{(k-1)d}$ Haar functions of level $k$. It is obvious that
$$
S_k^d = \mathrm{span}(\cup_{l=0}^k H_l^d),
$$
and that $H^d$ is a complete orthogonal system in $L_2(I^d)$. 

Since each Haar function in $H^d$ has support on a $d$-dimensional dyadic cube,
we sometimes call this system \emph{isotropic}, in contrast to the \emph{anisotropic} tensor-product Haar system $\tH^d$ defined in (\ref{HTPS}),
where the supports of the tensor-product Haar functions $\th\in \tH^d$ are $d$-dimensional dyadic rectangles. Note that $\tH^d$ can also be organized into blocks
$\tH_k^d$, where for $k\ge 1$ the block $\tH_k^d$ consists of the tensor-product Haar functions $\th\in S_k^d$ orthogonal to $S_{k-1}^d$, and $\tH_0^d=H_0^d=\{\chi_{I^d}\}$. Obviously, we have
$$
\mathrm{span}(\tH_k^d)=\mathrm{span}(H_k^d),\qquad k\in \mathbb{Z}_+,
$$
such that $\tH^d$ represents a level-wise transformation of $H^d$ vice versa.
The basis properties of $\tH^d$ in the spaces $B_{p,q,1}^s(I^d)$ are exhaustively dealt with 
in Section \ref{sec52}, see Theorem \ref{theo4}. It turns out that for $0<p< 1$ the two Haar systems $H^d$ and $\tH^d$ behave quite differently in this respect.

As for the univariate case, the ordering of the Haar functions within the blocks $H^d_k$ can be arbitrary. The statements of Theorems \ref{theo1} and \ref{theo2} hold for any enumeration of $H^d$ as long as the enumeration
does not violate the natural ordering by level $k$. Note that slightly more general orderings have been considered in \cite{GSU2019a} for Haar wavelet systems on $I^d$ and $\mathbb{R}^d$.

\subsection{Function spaces}\label{sec22}
The Besov function spaces $B_{p,q,1}^s(I^d)$ are traditionally defined for $s>0$, $0<p,q\le \infty$, as the set of all $f\in L_p(I^d)$ for which the quasi-norm
$$
\|f\|_{B_{p,q,1}^s}:=\left\{\ba{ll}(\|f\|_{L_p}^q+ \| t^{-s-1/q} \omega(t,f)_p\|_{L_q(I)}^q)^{1/q},& 0<q<\infty\\ \|f\|_{L_p}+ \sup_{t\in I} t^{-s} \omega(t,f)_p,&q=\infty  \ea\right.
$$
is finite. Formally, the definition makes sense for all $s\in \mathbb{R}$, see below for further comments in this direction. Here, 
$$
\omega(t,f)_p:=\sup_{0<|y|\le t} \|\Delta_y f\|_{L_p(I^d_y)}, \qquad t>0,
$$
stands for the first-order $L_p$ modulus of smoothness, and
$$
\Delta_y f(x):=f(x+y)-f(x), \qquad x\in I^d_y:=\{z\in I^d:\;z+y\in I^d\}, \quad y\in \mathbb{R}^d,
$$ 
denotes the first-order forward difference. Here and throughout the remainder of the paper, we adopt the following notational convention: If the domain is $I^d$,
we omit the domain in the notation for spaces and quasi-norms, e.g., we write
$B_{p,q,1}^s$ instead of $B_{p,q,1}^s(I^d)$, and $\|\cdot\|_{L_p}$ instead of
$\|\cdot\|_{L_p(I^d)}$. An exception are the formulation of theorems. Also, by $c,C$ we denote generic positive constants that may change from line to line, and, unless stated otherwise, depend on $p,q,s$ only. The notation $A\approx B$ is used if $cA\le B\le CA$ holds for two such constants $c,C$.

The space  ${B}_{p,q,1}^s$ is a quasi-Banach space equipped with
a $\gamma$-quasi-norm, where $\gamma=\min(p,q,1)$, meaning that $\|\cdot\|_{{B}_{p,q,1}^s}$
is homogeneous and satisfies
$$
\|f+g\|_{{B}_{p,q,1}^s}^{\gamma}\le \|f\|_{{B}_{p,q,1}^s}^{\gamma}+\|g\|_{{B}_{p,q,1}^s}^{\gamma}.
$$
Similarly, $L_p$ is a quasi-Banach space equipped with
a $\gamma$-quasi-norm if we set $\gamma=\gamma_p:=\min(p,1)$. 

If $0<q<\infty$ then the spaces ${B}_{p,q,1}^s$ are of interest only if $0\le s < 1/\gamma_p$.
Indeed, if $f\in  {B}_{p,q,1}^s$ for some $s\ge \gamma_p$, $0<q<\infty$, then using the properties of the first-order $L_p$ modulus of smoothness we have $\omega(t,f)_p=\mathrm{o}(t^{1/\gamma_p})$, $t\to 0$, which in turn implies $\omega(t,f)_p=0$ for all $t>0$ and $f(x)=\xi$ for some constant $\xi\in\mathbb{R}$ almost everywhere on $I^d$. Thus, in this case ${B}_{p,q,1}^s$ deteriorates to the set of constant functions on $I^d$. 
On the other hand, one has $\omega(t,f)_p\le 2^{1/\gamma_p}\|f\|_{L_p}$ which implies that 
$$
\|f\|_{{B}_{p,q,1}^s}\approx \|f\|_{L_p},\qquad f\in L_p,\qquad s<0.
$$
In other words, ${B}_{p,q,1}^s=L_p$ for all $s<0$ (this holds also for $s=0$ and $q=\infty$). 

To conclude this short discussion of the definition and properties of $B_{p,q,1}^s$-spaces, let us motivate our basic assumption  
(\ref{PR}) on the parameter range adopted in this paper. The case $s=0$ is in some sense exceptional and often not considered at all, we return to it in Section \ref{sec53}.
Since our main concern is the Schauder basis property of the countable Haar systems $H^d$ and $\tH^d$, we can also neglect all parameters for which ${B}_{p,q,1}^s$
is non-separable or does not contain piecewise constant functions on dyadic partitions. The separability requirement excludes the spaces with $p=\infty$ or $q=\infty$. Since, with constants also depending on $k$, we have
$$
\omega(t,h)_{p}\approx t^{1/p},\qquad t\to 0,
$$
for any Haar function $h\in H_k^d$, $k=1,2,\ldots$, we see that $H^d$ is not contained in $B_{p,q,1}^s$ whenever $s\ge 1/p$, $0<q<\infty$. Thus, the restrictions in (\ref{PR})  are natural.  

For completeness, we give the definition of the distributional Besov spaces using dyadic Fourier transform decompositions, see e.g. \cite[Section 1.1]{Tr2010}. Denote by $\mF: S'(\mathbb{R}^d) \to S'(\mathbb{R}^d)$ the Fourier transform operator on the set of tempered distributions. Consider a smooth partition of unity $\{\phi_k\}_{k\in \mathbb{Z}_+}$, where $\phi_0\in C^\infty(\mathbb{R}^d)$ satisfies $\phi_0(x)=1$ for $|x|\le 1$ and $\phi_0(x)=0$ for $|x|\ge 3/2$, and $\phi_k(x)=\phi_0(2^{-k}x)- \phi_0(2^{-k+1}x)$ for $x\in \mathbb{R}^d$,
and  $k\ge 1$. Then a tempered distribution
$f\in S'(\mathbb{R}^d)$ belongs to $B_{p,q}^s(\mathbb{R}^d)$ if the $\min(p,q,1)$-quasi-norm
$$%\be\label{DiBN}
\|f\|_{B_{p,q}^s(\mathbb{R}^d)}:=\| (\| 2^{ks}\mF^{-1}\phi_k\mF f\|_{L_p(\mathbb{R}^d)})_{k\in \mathbb{Z}_+}\|_{\ell_q(\mathbb{Z}_+)}
$$%\ee
is finite. By exchanging the order of taking $L_p(\mathbb{R}^d)$ and $\ell_q(\mathbb{Z}_+)$ quasi-norms
in (\ref{DiBN}), one defines the Triebel-Lizorkin spaces $F_{p,q}^s(\mathbb{R}^d)$.
Spaces on domains are defined by restriction. In particular, for the domain $I^d$
$$
B_{p,q}^s=\{f: \; \exists\,g\in B_{p,q}^s(\mathbb{R}^d) \mbox{ such that } f=g|_{I^d}\}
$$
and
$$
\|f\|_{B_{p,q}^s}:= \inf_{g:\;f=g|_{I^d}} \|g\|_{B_{p,q}^s(\mathbb{R}^d)}.
$$
This definition, and many equivalent ones, are surveyed in \cite{Tr2006} and \cite[Chapter 1]{Tr2010}, with references to earlier papers. In particular, the equivalence (\ref{PR31}) is mentioned in \cite[Section 1.1]{Tr2010}.

\subsection{Piecewise constant $L_p$-approximation}\label{sec23}
We start with introducing an equivalent quasi-norm in ${B}_{p,q,1}^s$ which is based on approximation techniques using piecewise constant approximation on dyadic partitions.  Let
$$
E_k(f)_p:= \inf_{s\in S_k^d} \|f-s\|_{L_p},\qquad k=0,1,\ldots,
$$
denote the best approximations to $f\in L_p$ with respect to $S_k^d$.
\begin{lem}\label{lem1} Let $0<p,q< \infty$, $0\le s<1/p$, and $d\ge 1$. Then
\be\label{Bnorm}
\|f\|_{{A}_{p,q,1}^s}:=\left(\|f\|_{L_p}^q + \sum_{k=0}^\infty (2^{ks} E_k(f)_p)^q\right)^{1/q}
\ee
provides an equivalent quasi-norm on ${B}_{p,q,1}^s$.
\end{lem}
This result follows from the direct and inverse inequalities relating best approximations
$E_k(f)_p$ and moduli of smoothness $\omega(t,f)_p$ which have many authors. In the univariate case  $d=1$, see e.g. Ul'yanov \cite{Ul1964}, Golubov \cite{Go1972} for $1\le p<\infty$, and \cite[Section 2]{SKO1975} for $0<p<1$.  Lemma \ref{lem1} is a partial case of  \cite[Theorem 5.1]{DP1988},  for $d=1$ and $0<p<1$ see \cite[Theorem 6]{Os1980}. The proofs for $s>0$ also cover the case $s=0$ not mentioned in these papers. Note that in \cite{DP1988} the parameter range $1\le s< 1/p$, $0<p<1$, is formally excluded but the result holds for the special case of piecewise constant approximation.  With the appropriate modification of the quasi-norm, such an approximation-theoretic characterization also holds
for $q=\infty$ and $0\le s<1/p$. 

The norm equivalence (\ref{Bnorm}) automatically implies that the set of all dyadic step functions 
$$%\be\label{Span}
S^d:=\mathrm{span}(H^d)=\mathrm{span}(\{S_k^d\}_{k\in \mathbb{Z}_+})
$$%\ee
is dense in ${B}_{p,q,1}^s(I^d)$ for the parameter values stated in Lemma \ref{lem1}. It can be used to prove sharp embedding theorems of ${B}_{p,q,1}^s$ into $L_r$.
In particular, we have continuous embeddings
\be\label{BtoL}
{B}_{p,p,1}^{d(1/p-1)} \subset {B}_{p,1,1}^{d(1/p-1)}\subset L_1,\qquad (d-1)/d < p < 1.
\ee
We refer to \cite{Os1980} for $d=1$, and to \cite[Theorem 7.4]{DP1988} for $d>1$. A local version of the associated embedding inequality will be used in Section \ref{sec3}.

At the heart of the counterexamples constructed in Section \ref{sec4} for $p\le 1$ is a simple observation about best $L_p$ approximation by constants
which we formulate as
\begin{lem}\label{lem2}
Let $(\Omega,\mathcal{A},\mu)$ be a finite measure space, and let the function $f\in L_p(\Omega):=L_p(\Omega,\mathcal{A},\mu)$, $0<p\le 1$, equal a constant $\xi_0$ on a measurable set $\Omega'\in \mathcal{A}$ of measure  $\mu(\Omega')\ge \frac12 \mu(\Omega)$. Then
$$
\|f-\xi_0\|_{L_p(\Omega)} = \inf_{\xi\in\mathbb{R}} \|f-\xi\|_{L_p(\Omega)},
$$
i.e., best approximation by constants in $L_p(\Omega)$ is achieved by setting $\xi=\xi_0$.
\end{lem}
{\bf Proof}. Indeed, under the above assumptions and by the inequality $|a+b|^p\le |a|^p+|b|^p$ we have
\bea
\|f-\xi\|_{L_p(\Omega)}^p&=&\int_{\Omega'}|f(x)-\xi|^p\,d\mu(x) +\mu(\Omega\backslash \Omega')|\xi-\xi_0|^p \\
&\ge& \int_{\Omega'}(|f(x)-\xi|^p +|\xi-\xi_0|^p)\,d\mu(x)\\
&\ge&
\int_{\Omega'}|f(x)-\xi_0|^p\,d\mu(x)=\|f-\xi_0\|_{L_p(\Omega)}^p
\eea
for any $\xi \in \mathbb{R}$, with equality for $\xi=\xi_0$. This gives the statement. \hfill $\Box$

\medskip
Note that the equivalence (up to constants depending on parameters but not on $f$) between $L_p$ quasi-norms and 
best approximations by constants holds also for $p\ge 1$ and under weaker assumptions on the relative measure of $\Omega'$
(e.g., $\mu(\Omega')/ \mu(\Omega)\ge \delta>0$ would suffice). We will apply this lemma to the Lebesgue measure on dyadic cubes in $I^d$ and special examples of dyadic step functions constructed below. Extensions to higher degree polynomial and spline approximation are possible as well (see the proof of the lemma on p. 535 in \cite{Os1981} for $d=1$).

\subsection{Schauder basis property}\label{sec24}
A sequence $(f_m)_{m\in \mathbb{N}}$ of elements of a quasi-Banach space $X$ is called a \emph{Schauder basis} in $X$ if every $f\in X$ possesses a 
\emph{unique} series representation
$$
f =\sum_{m=1}^\infty c_m f_m
$$
converging in $X$. If every rearrangement of $(f_m)_{m\in \mathbb{N}}$ is a Schauder basis in $X$ then this system is called \emph{unconditional Schauder basis}.
Below we will rely on the following criterion whose proof for Banach spaces can easily be extended to the quasi-Banach space case.
\begin{lem}\label{lem3} The sequence $(f_m)_{m\in \mathbb{N}}$ of elements of a quasi-Banach space $X$ is a Schauder basis in the quasi-Banach space $X$ if and only if its span is dense in $X$, and there exists a sequence $(\lambda_m)_{m\in \mathbb{N}}$ of continuous linear functionals on $X$ such that 
$$
\lambda_n(f_m)=\left\{\ba{ll} 1,&m=n,\\0,& m\neq n,\ea\right.\qquad m,n\in \mathbb{N},
$$
and the associated partial sum operators 
$$
S_n(x):= \sum_{m=1}^n \lambda_m(x) f_m,\qquad n\in \mathbb{N},\qquad x\in X,
$$
are uniformly bounded operators in $X$. \\
In order for $(f_m)_{m\in \mathbb{N}}$ to be unconditional, all operators
$$
S_J(x):= \sum_{m\in J} \lambda_m(x) f_m,\qquad x\in X,
$$
where $J$ is an arbitrary finite subset of $\mathbb{N}$, must be uniformly bounded operators in $X$.
\end{lem}
An immediate consequence of Lemma \ref{lem3} is that, in order to possess a Schauder basis at all, $X$ must have a sufficiently
rich dual space $X'$ of continuous linear functionals.

If $H^d$ is a Schauder basis in a quasi-Banach space $X$ of functions
defined on $I^d$ then  $S^d=\mathrm{span}(H^d)$ must be a dense subset 
of $X$ by the density condition in Lemma \ref{lem3}. This is satisfied for all $X={B}_{p,q,1}^{s}$ with parameters satisfying (\ref{PR}). Moreover, since 
$S^d\subset L_\infty\subset L_2$ and $H^d$ is an orthogonal system in $L_2$, any
dyadic step function $g\in S^d$ has a unique Haar expansion 
given by
\be\label{HE}
g= \sum_{h\in H^d} \lambda_h(g) h,\qquad \lambda_h(g):=2^{kd}\int_{I^d} g h \,dx.
\ee
Since for $g\in S^d$ only finitely many coefficients $\lambda_h(g)$
do not vanish, the summation in (\ref{HE}) is finite, and there are no convergence issues.
Thus, for the Schauder basis property of $H^d$ in $X$ to hold, the coefficient functionals $\lambda_h(g)$ in (\ref{HE}) must be extendable to elements in $X'$, and the level $k$ partial sum operators
\be\label{Pk}
P_kg = \sum_{l=0}^k \sum_{h\in H^d_l} \lambda_h(g) h, \qquad k=0,1,\ldots,
\ee
must form a sequence of \emph{uniformly} bounded linear operators in $X$. Due to the local support properties of the Haar functions within each block $H^d_l$ and the assumed ordering of the Haar wavelet system it often suffices to deal with this subsequence of partial sum operators. The statements in Theorem \ref{theo1} b) about the failure of the Schauder basis property
of $H^d$ in ${B}_{p,q,1}^{s}$ will be shown by either relying on Theorem \ref{theo2} or proving that the operators $P_k$ are not uniformly bounded.

Whenever $X$ is continuously embedded into $L_1$, the level $k$ partial sum operators $P_k$ extend to bounded projections with range $S_k^d$, and with constant values on the dyadic cubes in $T_k^d$ explicitly given by averaging. This comes in handy when computing $P_kf$ for concrete functions $f$. Indeed, the constant values taken by $P_kf$ on dyadic cubes in $T_k^d$ are given by 
\be\label{AvP}
P_kf(x)=av_\Delta(f):=2^{kd}\int_\Delta f\, dx, \qquad x\in \Delta,\quad \Delta\in T_k^d,,\quad k=0,1,\ldots,
\ee
if $f\in L_1(I^d)$. Note that the coefficient functionals $\lambda_h$ of the Haar expansion are finite linear combinations of functionals as defined in (\ref{AvP}), vice versa. Finally, for $X=L_2\subset L_1$ the level $k$ partial sum operator $P_k$ realizes the orthoprojection onto $S_k^d$.

\section{Proofs: Sufficient conditions} \label{sec3}
\subsection{The case $s=d(1/p-1)$} \label{sec31}
The positive results on the Schauder basis property of $H^d$ in $B_{p,q,1}^{d(1/p-1)}$  stated in Theorem \ref{theo1} b) for the parameter range 
\be\label{PR2}
\frac{d-1}{d} < p < 1,\quad 0<q\le p,
\ee
can be traced back to \cite{Os1981} for $d=1$, we 
reproduce the proof for $d\ge 1$ given in the preprint \cite{Os2018}.
For similar results in the case of distributional Besov spaces $B_{p,q}^s(\mathbb{R}^d)$ and $B_{p,q}^s$ we refer to \cite{Tr2010,Os2018,GSU2019a}.

According to (\ref{BtoL}), for the parameters in (\ref{PR2})
we have the continuous embedding
$$
{B}_{p,q,1}^{d(1/p-1)}\subset {B}_{p,p,1}^{d(1/p-1)}\subset L_1.
$$
This ensures that the Haar coefficient functionals $\lambda_h$ defined in (\ref{HE}) are continuous on ${B}_{p,q,1}^s$. Moreover, $S^d=\mathrm{span}(H^d)$ is dense in ${B}_{p,q,1}^s$.
Due to Lemma \ref{lem1} and Lemma \ref{lem3}, it is therefore sufficient to establish the inequality
\be\label{Unif}
\|Pg\|_{{A}_{p,q,1}^{d(1/p-1)}}\le C\|g\|_{{A}_{p,q,1}^{d(1/p-1)}},\qquad
g\in B_{p,q,1}^{d(1/p-1)},
\ee
for any partial sum operator $P$ of the Haar expansion (\ref{HE}), with a constant $C$ independent of $g$ and $P$, if the parameters satisfy
(\ref{PR2}).

According to our ordering convention for $H^d$, any partial sum operator $P$
can be written, for some $k=0,1,\ldots$ and some subset $\bar{H}^d_{k+1}\subset H^d_{k+1}$, in the form
\be\label{defPS}
{P}g=P_kg + \sum_{h\in \bar{H}^d_{k+1}} \lambda_h(g)h \in S_{k+1}^d.
\ee
For $\bar{H}^d_{k+1}=\emptyset$, we get $P=P_k$ as partial case.

The first step for establishing (\ref{Unif}) is the proof of the inequality
\be\label{Step11}
\|{P}g\|_{L_p}^p\le C 2^{kd(p-1)} \sum_{\Delta\in T_k^d} \|g\|_{L_1(\Delta)}^p,
\ee
with the explicit constant $C=2^d$.
By (\ref{AvP}), we have
$$
\|P_kg\|_{L_p}^p =\sum_{\Delta\in T_k^d} 2^{-kd} \left(2^{kd}\int_\Delta g \,dx\right)^p\le 2^{kd(p-1)} \sum_{\Delta\in T_k^d} \|g\|_{L_1(\Delta)}^p.
$$
The remaining $h\in \bar{H}^d_{k+1}$ can be grouped by their support cubes $\Delta\in T_k^d$. Each such group may hold up to $2^d-1$
Haar functions with the same $\mathrm{supp}(h)=\Delta\in T_k^d$. By the definition of the Haar coefficient functionals $\lambda_h(g)$ for each term $\lambda_h(g)h$ associated with such a group we obtain  the estimate
$$
\|\lambda_h(g)h\|_{L_p}^p = |\lambda_h(g)|^p \|h\|_{L_p(\Delta)}^p \le 2^{kdp}\|g\|_{L_1(\Delta)}^p\cdot 2^{-kd}=2^{kd(p-1)}\|g\|_{L_1(\Delta)}^p.
$$
Thus, using the $p$-quasi-norm triangle inequality for $L_p$, 
\bea
\|{P}g\|_{L_p}^p&\le& \|P_kg\|_{L_p}^p + \sum_{h\in \bar{H}^d_{k+1}} \|\lambda_h(g)h\|_{L_p}^p\\
&\le& 2^{kd(p-1)} \sum_{\Delta\in T_k^d} \|g\|_{L_1(\Delta)}^p +  \sum_{\Delta\in T_k^d} (2^d-1)2^{kd(p-1)}\|g\|_{L_1(\Delta)}^p\\
&=& 2^d 2^{kd(p-1)} \sum_{\Delta\in T_k^d} \|g\|_{L_1(\Delta)}^p,
\eea
we obtain (\ref{Step11}).

Now we apply the embedding inequality associated with (\ref{BtoL}), with the appropriate coordinate transformation, locally on $\Delta$ to the terms $\|g\|_{L_1(\Delta)}^p$. This gives
$$
\|g\|_{L_1(\Delta)}^p\le C\left(2^{kd(1-p)}\|g\|_{L_p(\Delta)}^p+\sum_{l=k}^\infty 2^{ld(1-p)}E_{l}(g)_{p,\Delta}^p\right)
$$
for each $\Delta\in T_k^d$, where 
$$
E_{l}(g)_{p,\Delta}:=\inf_{s\in S_{l}^d} \,\|g-s\|_{L_p(\Delta)},\qquad l=k,k+1,\ldots,
$$
denotes the local best $L_p$ approximation by dyadic step functions restricted to cubes $\Delta$ from $T_k^d$. Since
$$
\|g\|_{L_p}^p=\sum_{\Delta\in T_k^d} \|g\|_{L_p(\Delta)}^p,\qquad 
E_{l}(g)_{p}^p=\sum_{\Delta\in T_k^d} E_{l}(g)_{p,\Delta}^p,\quad l=k,k+1,\ldots,
$$
after substitution into (\ref{Step11}), we arrive at the estimate
\be\label{Step21}
\|{P}g\|_{L_p}^p\le C\left(\|g\|_{L_p}^p+2^{kd(p-1)}\sum_{l=k}^\infty 2^{ld(1-p)}E_{l}(g)_{p}^p\right)
\ee
for the $L_p$ quasi-norm of any partial sum $Pg$.

With the auxiliary estimate (\ref{Step21}) at hand, we turn now to the estimate of the 
$A_{p,q,1}^{d(1/p-1)}$ quasi-norm of $g-Pg$. Since $Pg\in  S_{k+1}^d$, we have
$$
E_{l}(g-Pg)_p=E_{l}(g)_p,\qquad l>k,
$$
while for $l\le k$ the trivial bound 
$$
E_{l}(g-Pg)_{p}\le \|g-Pg\|_{L_p(I^d)}
$$
will suffice. This gives
$$
\|g-Pg\|_{{A}_{p,q,1}^{d(1/p-1)}}^q=\|g-Pg\|_{L_p}^q + \sum_{l=0}^\infty (2^{ld(1/p-1)} E_l(g-Pg)_{p})^q\qquad\quad
$$
\be\label{Step31}
\qquad\qquad\qquad\le C\left(2^{kd(1/p-1)q}\|g-Pg\|_{L_p}^q + \sum_{l=k+1}^\infty (2^{ld(1/p-1)} E_l(g)_{p})^q\right),
\ee
uniformly for all $P$ and $g\in S^d$. Recall that $0<d(1/p-1)<1/p$ for the parameters in (\ref{PR2}).

To deal with the term $\|g-Pg\|_{L_p}$, we introduce the element $s_k\in S_k^d$ of best $L_p$ approximation, i.e.,
$$
\|g-s_k\|_{L_p}=E_k(g)_{p},
$$
and estimate with (\ref{Step21}) and $Ps_k=P_ks_k=s_k$ as follows:
\bea
\|g-Pg\|_{L_p}^p&\le& \|g-s_k\|_{L_p}^p+\|P(g-s_k)\|_{L_p}^p\\
&\le& \|g-s_k\|_{L_p}^p+C\left(\|g-s_k\|_{L_p(I^d)}^p+2^{kd(p-1)}\sum_{l=k+1}^\infty 2^{ld(1-p)}E_{l}(g-s_k)_{p}^p\right)
\eea
\be\label{Step3a}
\le C2^{kd(p-1)}\sum_{l=k}^\infty 2^{ld(1-p)}E_{l}(g)_{p}^p.\qquad\qquad\qquad\qquad\qquad
\ee

%Up to this point, only $g\in B_{p,p,1}^{sd(1/p-1)}\subset L_1$ has been assumed. 

Now we can finish the proof of (\ref{Unif}). According to (\ref{Step3a}) we get for the first term in the right-hand side of  (\ref{Step31})  
$$
2^{kd(1/p-1)q}\|g-Pg\|_{L_p}^q\le C\left(\sum_{l=k}^\infty 2^{ld(1-p)}E_{l}(g)_{p}^p
\right)^{q/p}\le C\sum_{l=k}^\infty (2^{ld(1/p-1)}E_{l}(f)_{p})^q,
$$
where the inequality
$$%\be\label{pinq}
\left(\sum_{l=0}^\infty a_l\right)^\gamma\le \sum_{l=0}^\infty a_l^\gamma, \qquad a_l\ge 0, \quad 0<\gamma \le 1,
$$%\ee
has been used with $\gamma=q/p\le 1$, $a_l=2^{ld(1-p)}E_{l}(g)_{p}^p$ for $l\ge k$, and $a_l=0$ for $l<k$.
After substitution into
(\ref{Step31}) we arrive at
\be\label{Unif1}
\|g-Pg\|_{{A}_{p,q,1}^{d(1/p-1)}}^q\le C\sum_{l=k}^\infty (2^{ls} E_l(g)_{p})^q\le C\|g\|_{{A}_{p,q,1}^{d(1/p-1)}}^q,\qquad g\in B_{p,q,1}^{d(1/p-1)},
\ee
for the parameters (\ref{PR2}).
Since the ${A}_{p,q,1}^{s}$-quasi-norm is a $q$-quasi-norm for the parameters in (\ref{PR2}),
(\ref{Unif1}) is equivalent with (\ref{Unif}). This proves the Schauder basis
property for $H^d$ in ${B}_{p,q,1}^{d(1/p-1)}$ for this  parameter range.

Under the condition (\ref{PR2}), unconditionality of $H^d$ does not hold, see Section  \ref{sec4} for a counterexample. The parameter range for which Theorem \ref{theo1} asserts that $H^d$ is an unconditional Schauder basis in ${B}_{p,q,1}^{s}$ is dealt with in the next subsection.
 
 \medskip
For use in the next subsection, we mention the following by-product of the considerations leading to
(\ref{Step3a}). Consider the $k$-th block 
$$
Q_kg:=(P_{k}-P_{k-1})g = \sum_{h\in {H}^d_{k}} \lambda_h(g)h, \qquad k=1,2,\ldots,
$$
of the Haar expansion (\ref{HE}). Using a standard compactness argument and the dyadic dilation- and shift-invariance of the Haar wavelet system $H^d$, we have
the equivalence of quasi-norms
$$
\|\sum_{h\in {H}^d_{k-1}: \mathrm{supp}(h)=\Delta} \gamma_h h\|_{L_p(\Delta)}^p 
\approx 2^{-kd}\sum_{h\in {H}^d_{k-1}: \mathrm{supp}(h)=\Delta} |\gamma_h|^p,\quad 0<p<\infty,
$$
which holds  with constants independent of the coefficient sequence $(\gamma_h)_{h\in H^d}$ and the  dyadic cubes $\Delta\in T_{k-1}^d$, $k=1,2,\ldots$.
Consequently,
\be\label{NE0}
\|\sum_{h\in {H}^d_{k-1}} \gamma_h h\|_{L_p}^p 
\approx 2^{-kd}\sum_{h\in {H}^d_{k-1}} |\gamma_h|^p,\quad 0<p<\infty,
\ee
with constants independent of $(\gamma_h)_{h\in H^d}$ and $k=1,2,\ldots$.
For the $L_p$ quasi-norm of $Q_kg$ this yields
$$
\|Q_kg\|_{L_p}^p=\|\sum_{h\in {H}^d_{k-1}} \lambda_h(g)h \|_{L_p}^p \approx
2^{-kd}\sum_{h\in {H}^d_{k-1}} |\lambda_h(g)|^p,\quad 0<p<\infty,
$$
Thus, we have
$$%\be\label{NE}
\sum_{h\in {H}^d_{k-1}} |\lambda_h(g)|^p \approx 2^{kd} \|Q_kg\|_{L_p}^p\le C2^{kd} (\|g-P_kg\|_{L_p}^p+\|g-P_{k-1}g\|_{L_p}^p),
$$%\ee
$k=1,2,\ldots$, for arbitrary $g\in L_1$ and $0<p<\infty$.

In combination with (\ref{Step3a}) this gives the estimate
\be\label{Step3b}
\sum_{h\in {H}^d_{k-1}} |\lambda_h(g)|^p \le C2^{kdp}\sum_{l=k-1}^\infty 2^{ld(1-p)}E_{l}(g)_{p}^p,\qquad k=1,2,\ldots,
\ee
for the $k$-th block of Haar coefficients of arbitrary $g\in B_{p,p,1}^{d(1/p-1)}$, $\frac{d-1}d < p < 1$. 

\subsection{Unconditionality}\label{sec32}
Now we turn to the parameter range
\be\label{PR3}
\max(d(1/p-1),0) < s < 1/p,\quad 0 < q < \infty,\quad (d-1)/d < p <\infty,
\ee
and establish the unconditional Schauder basis property of $H^d$ in ${B}_{p,q,1}^{s}$ by proving a slightly stronger statement. 
\begin{theo}\label{theo3}
	For the parameter range (\ref{PR3}), the mapping
	\be\label{LambdaMap}
	\Lambda :\; g \longmapsto \Lambda g := (\lambda_h(g))_{h\in H^d}, \qquad g\in L_1,
	\ee
	provides an isomorphism between ${B}_{p,q,1}^{s}(I^d)$ and a weighted $\ell_q(\ell_p)$ space, more precisely,
	we have
	\be\label{Iso}
	\|g\|_{{B}_{p,q,1}^{s}} \approx \|\Lambda g\|_{\ell_p(\ell_q)}:= \left(\sum_{k=0}^\infty \left(\sum_{h\in H^d_k} 2^{k(sp-d)}|\lambda_h(g)|^p\right)^{q/p}\right)^{1/q}.
	\ee
	Consequently, $H^d$ is an unconditional Schauder basis in ${B}_{p,q,1}^{s}(I^d)$ for the parameter range (\ref{PR3}).
\end{theo}
{\bf Proof}. For all parameters considered in (\ref{PR3}) we have the continuous embedding ${B}_{p,q,1}^{s}\subset L_1$. This ensures that $\Lambda$ is well-defined on ${B}_{p,q,1}^{s}$. In the following, we will concentrate on the case $(d-1)/d< p < 1$ 
which is partly new. For $1\le p<\infty$, the result is fully covered by \cite[Theorem 2.26 (i)]{Tr2010} stating the above isomorphism for $B_{p,q}^s$, since then (\ref{PR3}) implies that $B_{p,q}^s$ and $B_{p,q,1}^s$ coincide up to equivalent norms, see also \cite{Ro2016,Ro2016b}. Therefore, we only sketch the arguments that allow us to prove the result for $1\le p <\infty$ directly, without reference to the results for distributional Besov spaces $B_{p,q}^s$.

\emph{Step 1}. The upper bound
\be\label{Upper}
\|\Lambda g\|_{\ell_p(\ell_q)}\le C\|g\|_{{B}_{p,q,1}^{s}}, \qquad g\in {B}_{p,q,1}^{s},
\ee
can be proved as follows. For $(d-1)/d< p < 1$ we can use (\ref{Step3b}): Set temporarily 
$$
a_k:= \sum_{h\in {H}^d_{k}} |\lambda_h(g)|^p,\qquad k=0,1,\ldots .
$$
With this notation, we have
\bea
\|\Lambda g\|_{\ell_p(\ell_q)}^q&=&\sum_{k=0}^\infty (2^{k(sp-d)}a_k)^{q/p}\\
&\le& a_0^{q/p} + C \sum_{k=0}^\infty \left( 2^{k(sp-d(1-p))}\sum_{l=k}^\infty
2^{ld(1-p)}E_{l}(g)_{p}^p \right)^{q/p}
\eea
\be\label{Step4}
\qquad\qquad\qquad \le C\left(\|g\|_{L_p}^{q} + \sum_{k=0}^\infty \left( 2^{k(sp-d(1-p))}\sum_{l=k}^\infty
2^{ld(1-p)}E_{l}(g)_{p}^p \right)^{q/p}\right). 
\ee
In the last estimation step, we substituted the estimate
\bea
|a_0|^{q/p}&=&\left|\int_{I^d} g(x)\,dx \right|^q\le \|g\|_{L_1}^q \le C\|g\|_{A_{p,p,1}^{d(1/p-1)}}^q \\
& \le & C\left( \|g\|_{L_p}^p +\left(\sum_{l=0}^\infty 2^{ld(1-p)}E_{l}(g)_{p}^p\right)^{q/p}\right).
\eea
which follows from the definition of $a_0$ and the continuous embedding (\ref{BtoL}).

Now recall that under the assumption (\ref{PR3}) we have $d(1/p-1) < s< 1/p$. We can therefore choose $\epsilon$ such that $0 < \epsilon < s-d(1/p-1)$, and apply the Hardy-type inequality 
\be\label{Hardy}
\left(\sum_{l=k}^\infty b_l^r\right)^{1/r}\le C_{\epsilon,q/r}2^{-k\epsilon}\left(\sum_{l=k}^\infty (2^{l\epsilon}b_l)^q\right)^{1/q},\quad \epsilon >0, \quad k=0,1,\ldots,
\ee
valid for non-negative sequences $(b_l)_{l\in \mathbb{Z}_+}$ and all $0 <r,q <\infty$, with $b_l=2^{ld(1/p-1)}E_{l}(f)_{p}$ and $r=p$. This gives
$$
\left(\sum_{l=k} 2^{ld(1-p)}E_{l}(g)_{p}^p\right)^{q/p}\le C 2^{-k\epsilon q}\sum_{l=k}^\infty (2^{l(\epsilon+d(1/p-1))}E_{l}(g)_{p})^q, \quad k=0,1,\ldots.
$$
After substitution into (\ref{Step4}), we arrive at 
\bea
\|\Lambda g\|_{\ell_p(\ell_q)}^q&\le& 
C\left(\|g\|_{L_p}^{q} + \sum_{k=0}^\infty 2^{kq(s-d(1/p-1)-\epsilon)}\sum_{l=k}^\infty (2^{l(\epsilon+d(1/p-1))}E_{l}(g)_{p})^q\right)\\
&=& C\left(\|g\|_{L_p}^{q} + \sum_{l=0}^\infty (2^{l(\epsilon+d(1/p-1))}E_{l}(g)_{p})^q \sum_{l=0}^k 2^{kq(s-d(1/p-1)-\epsilon)} \right)\\
&\le& C\left(\|g\|_{L_p}^{q} + \sum_{l=0}^\infty (2^{ls}E_{l}(g)_{p})^q \sum_{k=0}^l 2^{kq(s-d(1/p-1)-\epsilon)} \right) \le C\|g\|_{A_{p,q,1}^s}^q.
\eea
This proves (\ref{Upper})  for the range $(d-1)/d <p<1$. 

For $1\le p<\infty$ we can use the inequality
$$
\|g-P_k g\|_{L_p} \le C E_k(g)_p,\qquad g\in L_p,\qquad k\in \mathbb{Z}_+,
$$
which follows from the uniform boundedness of the projectors $P_k$ in $L_p$
(for details, see Section \ref{sec53}).
As above, in conjunction with (\ref{NE}) this gives 
$$
a_k=\sum_{h\in {H}^d_{k}} |\lambda_h(g)|^p\le C2^{kd}\|Q_kg\|_{L_p}^p\le C2^{kd}E_{k-1}(g)_p^p,\qquad k=1,2,\ldots,
$$
and $a_0\le \|g\|_{L_p}^p$. Then (\ref{Upper}) follows by simple substitution into the expression for the weighted $\ell_q(\ell_p)$-norm of $\Lambda g$ defined  in (\ref{Iso}).

It remains to prove that $\Lambda$ is surjective. Since the operator $\Lambda$ is obviously injective on any Besov space $B_{p,q,1}^s$ embedded into $L_1$, surjectivity together with (\ref{Upper}) automatically implies boundedness of the inverse mapping $\Lambda^{-1}$
as a consequence of the open mapping theorem for $F$-spaces (all $\gamma$-quasi-normed
Banach spaces and, in particular, the Besov spaces $B_{p,q,1}^s$ are $F$-spaces).

Let $\Gamma=(\gamma_h)_{h\in H^d}$ be an arbitrary sequence with finite weighted $\ell_q(\ell_p)$-quasi-norm:
\be\label{Sequ}
\|\Gamma\|_{\ell_q(\ell_p)}:= \left(\sum_{k=0}^\infty \left(\sum_{h\in H^d_k} 2^{k(sp-d)}|\gamma_h|^p\right)^{q/p}\right)^{1/q} <\infty.
\ee
To show the surjectivity of $\Lambda$, we need to find a $g\in B_{p,q,1}^s$ such that
\be\label{Coff}
\lambda_h(g)=\gamma_h, \qquad h\in H^d.
\ee
We will show this in all detail only for the case $(d-1)/d < p < 1$ in (\ref{PR3}),
the case $1\le p<\infty$ is analogous, see \cite{Ro1976} for $d=1$ and \cite{Tr2010} for $d>1$.

Set
$$
q_k:=\sum_{h\in H_k^d} \gamma_h h,\quad p_k:= \sum_{l=0}^k q_l,\qquad k\in \mathbb{Z}_+.
$$
We first show that $(p_k)_{k\in \mathbb{Z}_+}$ is fundamental in $L_1$ whenever $s>d(1/p-1)$ in (\ref{Sequ}). Indeed, by (\ref{NE0}) and (\ref{Hardy}) we have
\bea
\|p_m-p_k\|_{L_1}&\le& \sum_{l=k+1}^\infty \|q_l\|_{L_1}\le C\sum_{l=k+1}^m 2^{-ld}\sum_{h\in {H}^d_{l-1}} |\gamma_h|%\\
\le C\sum_{l=k}^\infty 2^{-ld}\left(\sum_{h\in {H}^d_{l}} |\gamma_h|^p\right)^{1/p}\\
&\le& C2^{-k(s-d(1/p-1))}\left(\sum_{l=k}^\infty \left(\sum_{h\in H^d_{l}} (2^{l(sp-d)}|\gamma_h|^p\right)^{q/p}\right)^{1/q}\\
&\le& C2^{-k(s-d(1/p-1))}\|\Gamma\|_{\ell_q(\ell_p)}
\eea
for arbitrary $1\le k < m<\infty$. The Hardy-type inequality (\ref{Hardy}) has been applied in the last but one estimation step
with $\epsilon=s-d(1/p-1)>0$ to the sequence
$$
b_l=2^{-ld}\left(\sum_{h\in {H}^d_{l}} |\gamma_h|^p\right)^{1/p},\qquad l\ge k,
$$
whereas the parameter $r$ in (\ref{Hardy}) has been set to $r=1$. The above bound on 
$\|p_m-p_k\|_{L_1}$ shows that $(p_k)_{k\in\mathbb{Z}_+}$ is fundamental in $L_1$, and that it converges to a function $g\in L_1$. Obviously, this implies the $L_1$ convergence of the Haar series with coefficients $\gamma_h$ to $g$, i.e., 
$$
g=\sum_{h\in H^d} \gamma_h h,
$$
assuming the agreed upon blockwise ordering of the Haar functions.
By orthogonality of $H^d$ and $\lambda_h\in L_\infty=L_1'$ for all Haar functions, we also must have (\ref{Coff}). 

It remains to show that the above $g$ belongs to $B_{p,q,1}^s$. We use 
$$
E_k(g)_p^p\le \|g-p_k\|_{L_p}^p\le \sum_{l>k}\|q_l\|_{L_p}^p, \qquad k\in\mathbb{Z}_+,
$$
and, similarly, 
$$
\|g\|_{L_p}^p\le \sum_{l=0}^\infty\|q_l\|_{L_p}^p.
$$
As above, together with Lemma \ref{lem1}, (\ref{NE0}), and (\ref{Hardy}) with $0<\epsilon <s$, this gives
\bea
\|g\|_{B_{p,q,1}^s}^q&\le& C\|g\|_{A_{p,q,1}^s}^q\le C\sum_{k=0}^\infty \left(2^{ksp}\sum_{l=k}^\infty \|q_l\|_{L_p}^p\right)^{q/p} \\
&\le& C\sum_{k=0}^\infty 2^{k(s-\epsilon)q}\left(\sum_{l=k}^\infty 2^{l\epsilon q}\|q_l\|_{L_p}^q\right) \le C \sum_{l=0}^\infty\|q_l\|_{L_p}^q \\
&\le& C \sum_{l=0}^\infty  2^{lsq}\left(2^{-ld}\sum_{h\in {H}^d_{l}} |\gamma_h|^p\right)^{q/p}
=C\|\Gamma\|_{\ell_q(\ell_p)}^q.
\eea
This finishes the proof of Theorem \ref{theo3} for $(d-1)/d<p< 1$. 

The proof of the surjectivity of the map $\Lambda$ for $1\le p <\infty$ is similar.
We mention the minor differences.
Instead of $L_1$ convergence, we can directly prove $L_p$ convergence of the Haar series
with coefficient sequence $\Gamma$ since
\bea
\|p_m-p_k\|_{L_p}&\le& \sum_{l=k+1}^\infty \|q_l\|_{L_p}\le C2^{-ks}
\left(\sum_{l=k+1}^\infty (2^{ls} \|q_l\|_{L_p})^q\right)^{1/q}\\
&\le& C2^{-ks}\sum_{l=k}^\infty \left(2^{l(sp-d)}\sum_{h\in {H}^d_{l-1}} |\gamma_h|^p\right)^{q/p}, \qquad 1\le k<m <\infty,
\eea
where (\ref{Hardy}) has been applied with $\epsilon=s>0$, $b_l=\|q_l\|_{L_p}$, and $r=1$, followed by (\ref{NE0}). In the proof of $g\in B_{p,q,1}^s$ instead of the $p$-quasi-norm property for $p<1$, we use the additivity of the norm for $p\ge 1$:
$$
E_k(g)_p\le \sum_{l>k} \|q_l\|_{L_p},\qquad \|g\|_{L_p}\le  \sum_{l=0}^\infty \|q_l\|_{L_p}.
$$ 
The rest of the proof is, up to obvious changes in the application of (\ref{Hardy}), the same
as for $p <1$.

\section{Proofs: Necessary conditions}\label{sec4}
\subsection{Proof of Theorem \ref{theo2}} \label{sec41}
Let the parameters satisfy (\ref{PR1}). For $(d-1)/d<p<1$ the Besov space $B_{p,q,1}^s$ is continuously embedded into $L_1$ if and only if
$$%\be\label{PR4}
d(1/p-1)<s<1/p,\; 0<q<\infty\quad\mbox{or} \quad s=d(1/p-1), \;q\le 1, 
$$%\ee
see (\ref{BtoL}), while for $p\ge 1$ this embedding is obvious. In all these cases, the dual of $B_{p,q,1}^s$ 
is therefore infinite-dimensional. This proves the necessity of the condition in Theorem 
\ref{theo2}.
 
For the sufficiency, we can  concentrate on the range $0<p< 1$, and assume on the contrary that there is a non-trivial continuous linear functional $\phi$ on
$B_{p,q,1}^s$. Since the span of characteristic functions $\chi_\Delta$ of all dyadic 
cubes in $I^d$ is dense in $B_{p,q,1}^s$, there must be a dyadic cube $\Delta_0$ such that
$\phi(\chi_{\Delta_0})\neq 0$. Without loss of generality, we may assume that
$$
\Delta_0=I^d,\qquad \phi(\chi_{\Delta_0})=1.
$$
Each dyadic cube in $T_k^d$, $k\ge 0$, is the disjoint union of $2^d$ dyadic cubes from
$T_{k+1}^d$. Therefore, using the  linearity of $\phi$, we can construct by induction a sequence of dyadic cubes $\Delta_k\in T_k^d$ such that
$$%\be\label{DQ}
\Delta_0\supset \Delta_1\supset\ldots\supset \Delta_k\supset \ldots,\qquad \phi(\chi_{\Delta_k})\ge 2^{-kd},\quad k=0,1,\ldots.
$$%\ee

For a given sequence $a=(a_k)_{k\in \mathbb{Z}_+}$, consider the function
\be\label{Fm}
f_m = \sum_{l=0}^m a_l \chi_{\Delta_l} \in S^d_m, \qquad m=0,1,\ldots.
\ee
The $L_p$ and Besov space quasi-norms of $f_m$ can be computed exactly. Indeed, $f_m$ is constant on all cubes $\Delta\in T_k^d$ but $\Delta_k$, and 
equals the constant $\xi_k:=\sum_{l=0}^k a_l$ on $\Delta'_k:=\Delta_k\backslash \Delta_{k+1}$,
where 
$$
\mu(\Delta'_k)=(1-2^{-d})\mu(\Delta_k)=(1-2^{-d})2^{-kd}\ge \frac12 \mu(\Delta_k),
$$
and $\mu$ denotes the Lebesgue measure on $\mathbb{R}^d$.
Thus,   
\be\label{FL}
\|f_m\|_{L_p}^p = (1-2^{-d})\sum_{n=0}^m 2^{-nd}\left|\sum_{l=0}^n a_l\right|^p.
\ee
Moreover, by Lemma \ref{lem2} applied to $\Omega=\Delta_k$ and $\Omega'=\Delta'_k$ we have
\be\label{FE}
E_k(f_m)_p^p=\|f_m-\xi_k\|_{L_p(\Delta_k)}^p = (1-2^{-d})\sum_{n=k+1}^m 2^{-nd}\left|\sum_{l=k+1}^n a_l\right|^p,\quad k=0,1,\ldots,m-1,
\ee
while obviously $E_k(f_m)_p=0$ for $k\ge m$.

Now we choose the particular sequence 
$$
a_k=2^{kd}(k+1)^{-1}, \qquad k=0,1,\ldots,
$$
in (\ref{Fm}). Substituting into (\ref{FL}) and (\ref{FE}), we have for $0<p<1$
$$
\|f_m\|_{L_p}^p\approx \sum_{n=0}^m 2^{-nd}\left(\sum_{l=0}^n \frac{2^{ld}}{l+1}\right)^p
	\approx \sum_{n=0}^m \frac{2^{-nd(1-p)}}{(n+1)^p}\approx 1,
$$
and similarly
$$
E_k(f_m)_p^p\approx \sum_{n=k+1}^m \frac{2^{-nd(1-p)}}{(n+1)^p} \approx \frac{2^{-kd(1-p)}}{(k+1)^p},\qquad k=0,1,\ldots,m-1.
$$
Note that by the same token 
$$
\|f_m\|_{L_1}\approx \sum_{n=0}^m \frac{1}{n+1}\approx \ln(m+1),\qquad m=0,1,\ldots,
$$
and
\be\label{Phi}
\phi(f_m)=\sum_{l=0}^m a_l \phi(\chi_{\Delta_l})\ge \sum_{l=0}^m \frac{2^{ld}}{l+1}2^{-ld}
\ge c\ln(m+1) \to \infty,\qquad m\to \infty.
\ee

Next we compute the Besov quasi-norm of $f_m$:
\bea
\|f_m\|_{B_{p,q,1}^{s}}^q &\le& C\|f_m\|_{B_{p,q,1}^{s}}^q = C\left(\|f_m\|_{L_p}^q+\sum_{k=0}^{m-1} (2^{ks}E_k(f_m)_p)^q\right)\\
&\le& C\left(1+\sum_{k=0}^{m-1} \left(\frac{2^{k(s-d(1/p-1))}}{k+1}\right)^q\right)\\
&\le& C\sum_{k=0}^{\infty} \frac{2^{k(s-d(1/p-1))q}}{(k+1)^q}.
\eea
Thus, if $0<s<d(1/p-1)$ and $0<q<\infty$ or if $s=d(1/p-1)$ and $1<q<\infty$,
the sequence $(f_m)_{m\in \mathbb{Z}_+}$ is uniformly bounded in $B_{p,q,1}^s$
which according to (\ref{Phi}) contradicts the assumed boundedness of the linear functional of $\phi$. This proves the necessity part of Theorem \ref{theo2}.
We note that for $s<d(1/p-1)$ simpler examples can be used 
to show the same result, see \cite{Os1981,Os2018}.

\subsection{The case $s=d(1/p-1)>0$, $0<q\le p<1$}\label{sec42}
In Section \ref{sec31} we established that $H^d$ has the Schauder basis property for $B_{p,q,1}^{d(1/p-1)}$ if $0<q\le p$, $(d-1)/d<p<1$. That $H^d$ is not an unconditional Schauder basis for these spaces can be shown by a simple example. It also appears in \cite[Section 13]{GSU2019a} in similar context, a related construction for $d=1$ can be 
found in \cite{Kr1982}.

For $(d-1)/d <p<1$, consider the sequence
$$
f_m=2^{md}\chi_{\Delta_m}\in S_m^d,  \qquad \Delta_m:=[0,2^{-m})^d\in T_m^d,\quad m=0,1,\ldots.
$$
Obviously, $E_k(f_m)_p=0$ for $k\ge m$, and using Lemma \ref{lem2} we compute
$$
E_k(f_m)_p=\|f_m\|_{L_p} =2^{-md/p}2^{md}= 2^{-md(1/p-1)}, \quad k=0,1,\ldots,m-1.
$$
Thus, since
\be\label{FBN}
\|f_m\|_{B_{p,q,1}^{d(1/p-1)}}^q \approx \|f_m\|_{A_{p,q,1}^{d(1/p-1)}}^q=\|f_m\|_{L_p}^q\left(1 +\sum_{k=0}^{m-1} 2^{kd(1/p-1)q}\right)\approx 1,\qquad m=0,1,\ldots,
\ee
the $B_{p,q,1}^{d(1/p-1)}$ quasi-norms of $f_m$ are uniformly bounded for the indicated
parameter range.

The Haar coefficients of $f_m$ are easily computed:
$$
\lambda_h(f_m)=\left\{\ba{ll} 1,& h\in H_0^d,\\  2^{(k-1)d},& \Delta_m\subset \mathrm{supp}(h),\;h\in H_k^d,\, k=1,\ldots,m\\ 0,& \mbox{otherwise}. \ea\right.
$$
From this, it is immediate that 
$$
\sum_{h\in H^d_l} \lambda_h(f_m) h =f_l-f_{l-1},\qquad l=0,1,\ldots,m,
$$
where we have set $f_{-1}=0$ for convenience. Consider now the function
$$%\be\label{Gk}
g_{2k} :=\sum_{l=0}^k \sum_{h\in H^d_{2l}} \lambda_h(f_m) h = \sum_{l=0}^{2k} (-1)^lf_l 
= \sum_{l=0}^{2k} (-1)^l2^{ld}\chi_{\Delta_l},
$$%\ee
where $2k\le m$. This function is of the same type as the functions $f_{2k}$ considered in the previous subsection but with a different coefficient sequence $a=((-1)^l2^{ld})_{l\in \mathbb{Z}_+}$ and a specific nested sequence of dyadic cubes $\Delta_l\in T_l^d$. Using the formula (\ref{FE}), we compute
$$
E_l(g_{2k})_p\approx 2^{-ld(1/p-1)},\quad l=0,\ldots,2k-1, \qquad E_l(g_{2k})_p=0,\quad l\ge k,
$$
and conclude that
$$
\|g_{2k}\|_{B_{p,q,1}^{d(1/p-1)}}^q\ge c \sum_{l=0}^{2k-1} (2^{ld(1/p-1)}E_l(g_{2k})_p)^q\ge
ck,\qquad 2k \le m,
$$
grows unboundedly if $k,m\to \infty$, independently of $0<q<\infty$ and, in particular, for the range $0<q\le p$ of interest. In conjunction with (\ref{FBN}) this contradicts the unconditionality of $H^d$ since the $g_{2k}$ are partial sums of the Haar expansion of $f_m$ with respect to a specific finite subset
of $H^d$. Recall that by Lemma \ref{lem3} unconditionality of a Schauder basis requires the \emph{uniform} boundedness of the partial sum operators for arbitrary finite subsets  of basis elements.

\subsection{The case $s=d(1/p-1)>0$, $(d-1)/d < p< q \le 1$}\label{sec43}
For the parameters
\be\label{PR5}
s=d(1/p-1), \qquad (d-1)/d < p< q \le 1,
\ee
the failure of the Schauder basis property claimed in Theorem \ref{theo1} b) cannot
be deduced from Theorem \ref{theo2}. Since the Haar coefficient functionals  $\lambda_h$ are continuous on $B_{p,q,1}^s$, the finite-rank partial sum operators $P_k$ defined in (\ref{Pk}) represent bounded operators. However, they are not uniformly bounded on $B_{p,q,1}^s$ for the parameter values in (\ref{PR5}).

To establish this fact, we need another type of examples which were introduced in \cite{Os2018} in a slightly different form. Fix a $k\ge 1$ arbitrarily,
and select a subset $T'\subset T_k^d$ such that  each dyadic cube $\Delta\in T_{k-1} ^d$ contains  exactly $2^{d-1}$ cubes $\Delta'\in T'$ (and thus exactly $2^{d-1}$ cubes from $T^d_k\backslash T'$). This way  $T'$ contains $2^{kd-1}$ cubes from $T_k^d$ which will be enumerated in an arbitrary order, and denoted by $\Delta'_i $, $i=1,\ldots,2^{kd-1}$.
For each $i=1,\ldots,2^{kd-1}$, choose a dyadic cube $\Delta_i\in T_{k+i}^d$ contained in 
$\Delta'_i$, and set
$$%\be\label{Fk1}
f_k:= \sum_{i=1}^{2^{kd-1}} b_{i}\chi_{\Delta_i},\qquad b_i:=2^{(k+i)d}i^{-\alpha},\quad
i=1,\ldots,2^{kd-1},
$$%\ee
where $\alpha < 1/q$ is fixed.  

It is clear by this construction that on each dyadic cube $\Delta$ we have either $f_k=0$ on a subset $\Omega'\subset \Delta$ of measure $\mu(\Omega')\ge \frac12 \mu(\Delta)$, or $f_k$ is constant on $\Delta$. For $\Delta\in T_l^d$ with $l<k$, only the former option is possible due to the properties of $T'$. For $\Delta\in T_l^d$ with $l\ge k$, we have three mutually exclusive cases: 1) $\Delta\subset \Delta_i$ for some $i$, 2) $\Delta_i$ is strictly contained in $\Delta$ for some $i$, or 3) $\Delta$ does not intersect any of the $\Delta_i$. In case 1) and 3) $f_k$ is constant on $\Delta$, while in case 2) $f_k=0$ on a set $\Omega'\subset \Delta$ of
measure $\mu(\Omega')\ge (1-2^{-d})\mu(\Delta) \ge \frac12 \mu(\Delta)$. This enables 
the use of Lemma \ref{lem2} for the computation of best approximations. Since $p<1$, for $l=0,\ldots,k$, we obtain
$$
E_l(f_k)_p^p=\|f_k\|_{L_p}^p= \sum_{i=1}^{2^{kd-1}} 2^{-(k+i)d}b_{i}^p =\sum_{i=1}^{2^{kd-1}}
2^{-(k+i)d(1-p)}i^{-\alpha p}\approx 2^{-kd(1-p)}.
$$
For $l=k+1,\ldots,k+2^{kd-1}-1$ we similarly have
$$
E_l(f_k)_p^p=\|\sum_{i=l+1-k}^{2^{kd-1}} b_{i}\chi_{\Delta_i}\|_{L_p}^p=
\sum_{i=l+1-k}^{2^{kd-1}} 2^{-(k+i)d(1-p)}i^{-\alpha p}\approx  2^{-ld(1-p)}(l-k)^{-\alpha p},
$$
while $E_l(f_k)_p=0$ for $l\ge k+2^{kd-1}$. Substitution into the $A_{p,q,1}^{d(1/p-1)}$ quasi-norm expression gives
\bea
\|f_k\|_{B_{p,q,1}^{d(1/p-1)}}^q &\approx& 2^{-kd(1/p-1)q}\sum_{l=0}^k 2^{ld(1/p-1)q} +
\sum_{l=k+1}^{k+2^{kd-1}-1} (2^{ld(1/p-1)} 2^{-ld(1/p-1)} (l-k)^{-\alpha})^q\\
&\approx& 1 + \sum_{i=1}^{2^{kd-1}} i^{-\alpha q}\approx2^{kd(1-\alpha q)},\qquad k=1,2,\ldots,
\eea
where $\alpha q < 1$ has been used.

The construction of $T'$ also allows us to compute the best approximations of the partial sum $P_kf_k$ of the finite Haar expansion of $f_k$.
Indeed, according to (\ref{AvP}-\ref{Av}), $P_kf_k$ is given by
$$
P_kf_k(x)=\left\{\ba{ll} 2^{kd}2^{-(k+i)d}b_i=2^{kd} i^{-\alpha},& x\in \Delta'_i, \;i=1,\ldots,2^{kd-1},\\ 0, & x\in \Delta',\, \Delta'\in T_k^d\backslash T'.\ea \right.
$$
Since $P_kf_k\in S_k^d$ we have $E_l(P_kf_k)_p=0$ for $l\ge k$. By the selection rule of 
the cubes $\Delta'_i\in T'$, on any cube $\Delta\in T_l^d$ with $l<k$ we have $P_kf_k=0$ on the union of all $\Delta'\in T^d_k\backslash T'$ intersecting with $\Delta$ which has always exactly half of the measure of $\Delta$. For this reason, Lemma \ref{lem2} is applicable and gives
$$
E_l(P_kf_k)_p^p=\|P_kf_k\|_p^p= \sum_{i=1}^{2^{kd-1}} 2^{-kd}(2^{kd}i^{-\alpha})^p = 
2^{-kd(1-p)}\sum_{i=1}^{2^{kd-1}} i^{-\alpha p}
$$
for $l=0,\ldots,k-1$. Consequently, we have
$$
\|P_kf_k\|_{B_{p,q,1}^{d(1/p-1)}}^q \approx 2^{kd(1/p-1)q}\|P_kf_k\|_p^q \approx \left(\sum_{i=1}^{2^{kd-1}} i^{-\alpha p}\right)^{q/p}
\approx 2^{kd(1/p-1)q},
$$
since for the parameter range of interest $\alpha <1/q <1/p$.

Comparing the above estimates for the Besov quasi-norms of $P_kf_k$ and $f_k$ shows that
$$%\be\label{Unb}
\|P_k\|_{B_{p,q,1}^{d(1/p-1)}\to B_{p,q,1}^{d(1/p-1)}} \ge \frac{\|P_kf_k\|_{B_{p,q,1}^{d(1/p-1)}}}{\|f_k\|_{B_{p,q,1}^{d(1/p-1)}}}\ge c2^{kd(1/p-1/q)},\qquad k\ge 1,
$$%\ee
for values $0< p < q$. I.e., the partial sum operators $P_k$ are not uniformly bounded, and  $H^d$ cannot be a Schauder basis in $B_{p,q,1}^{d(1/p-1)}$ for the parameter range (\ref{PR5}). This concludes the proof of the necessity of the conditions in Theorem \ref{theo1}.

\section{Remarks and extensions}\label{sec5}

\subsection{Higher-order spline systems}\label{sec51} 
The Schauder basis property of spline systems 
of order $m>1$ in Lebesgue-Sobolev spaces and in Besov spaces $B_{p,q,m}^s$ defined by $m$-th order moduli of smoothness
$$
\omega_m(t,f)_p:=\sup_{0\le |y|\le t} \|\Delta_y^m f\|_{L_p(I^d_{y,m})}, \qquad t>0,
$$
where $\Delta^m_y$ denotes the $m$-th order forward difference operator with step-size $y$, and $I^d_{y,m} =\{x\in I^d:\; x+my\in I^d\}$,
has also attracted attention, see \cite{CiDo1972,CiFi1982,CiFi1983a,CiFi1983b} for $1\le p \le \infty$. Note that, with the exception of \cite{CiFi1983b}, in these papers the case $d>1$ has been treated by the tensor product construction, i.e., with analogs of $\tH^d$, not the Haar wavelet system $H^d$. For $0<p<1$, we refer to \cite{Os1981} which treats orthogonal spline systems in the one-dimensional case. Higher-order spline wavelet systems as unconditional Schauder bases in Besov-Triebel-Lizorkin spaces $B_{p,q}^s(\mathbb{R}^d)$ and $F_{p,q}^s(\mathbb{R}^d)$ have been considered in \cite[Section 2.5]{Tr2010}. 

Without going into much detail, we claim that analogs of the above results hold for $B_{p,q,m}^s$ and $m$-th order orthogonal spline wavelet systems in the parameter range
$$
0< s <\min(m,m-1+1/p),\qquad 0<p,q <\infty,
$$
and can be proved following the above reasoning for the Haar system $H^d$.
To be more specific, take the univariate orthogonal Ciesielski system $F^m$
of order $m\ge 1$. This is the system $\{f_n^{m-2},\,n\ge -m+2\}$ from \cite{CiDo1972,Ci1975}
obtained by Gram-Schmidt orthogonalization of a suitably chosen system of B-splines
of degree $m-1$ associated with the dyadic partitions $T_k$ of the unit interval $I$, $k=0,1,\ldots$. Its $d$-dimensional wavelet counterpart $F^{m,d}$ is then constructed along the lines of \cite[Section 10]{CiFi1983b} or \cite[Section 2.5.1]{Tr2010}, where the $f_n^{m}$ play the role of the wavelet functions, and the B-splines the role of the scaling functions, respectively. Note that what we call order $m$ in this paper means degree
$m-1$ and  $C^{m-2}$ smoothness  of the splines, respectively, and is used differently in the cited papers. Then, using the exponential decay properties of the Ciesielski functions $f_n^{m-2}$,
and the characterization of $B_{p,q,m}^s$ in terms of best approximations from \cite{DP1988} valid for the indicated range of smoothness parameters $s$, one first proves an analog 
$$
\|{P}^mg\|_{L_p}^p\le C 2^{kd(p-1)} \sum_{\Delta\in T_k^d} \|g\|_{L_1(\Delta)}^p,\qquad 0<p\le 1,\qquad g\in L_1,
$$
of the crucial estimate (\ref{Step11}), now for partial sum operators $P^m$ of level $k$ with respect to $F^{m,d}$, and then follows the proof of Theorem 1 in Section \ref{sec31}. This essentially leads to positive results for $\max(0,d(1/p-1)) < s <\min(m,m-1+1/p)$ and $0< q<\infty$ as well as on the critical line $s=d(1/p-1)$, $0<q\le p < 1$, for $d=1$ see
\cite{Os1981}.

The counterexamples for the remaining cases in the range $0<s\le d(1/p-1)$, $0<p<1$, can be built using linear combinations of B-splines with well-separated supports, as done in \cite{Os1981} for $d=1$ (see the lemma on p. 535 there). Some additional technical difficulties arise from the fact that partial sum operators are not as local as in the Haar case ($m=1$) but can be overcome using the exponential decay of the associated operator kernels in conjunction with the support separation in the examples.

Since for $m\ge 2$ and $d>1$ we have
$$
F^{m,d}\not\subset B_{p,q,m}^s,\qquad s\ge m, \quad 0<q<\infty,\quad 0<p<1,
$$
the restriction to the range $0<s<m$ for $0<p<1$ is natural. This is in contrast to the case $m=1$, where
$H^d$ remained an unconditional Schauder basis also for the range $1\le s <1/p$, $0<p<1$.
One could therefore ask if the spaces $B_{p,q,m}^s$ permit specific spline bases
for the remaining values 
\be\label{PR6}
m\le s <m-1+1/p, \qquad 0<q<\infty, \qquad 0<p<1, 
\ee
also for $d>1$. In order for a spline system
to belong to $B_{p,q,m}^s$ for this parameter range, it is desirable that it consists
of splines which are locally polynomials of exactly \emph{total} degree $m-1$, and globally
belong to $C^{m-2}$ over well-shaped and refinable partitions. These are the spline functions that are maximally smooth in $L_p$, $0<p<1$,
in the sense that their $m$-th order modulus of smoothness decays at the best possible rate $\mathrm{O}(t^{m-1+1/p})$, $t\to 0$. Tensor-product constructions such as $F^{m,d}$ lead to
$C^{m-2}$ smooth splines with local \emph{coordinate} degree $m-1$ but total degree
$d(m-1)>m-1$ for which only a $\mathrm{O}(t^{m})$ decay of the $m$-th order modulus of smoothness can be expected. To the best of our knowledge, spline systems with all properties desirable for the  construction of
Schauder bases in $B_{p,q,m}^s$ with parameters satisfying (\ref{PR6}) are available only in special cases. E.g., for $m=2$ semi-orthogonal prewavelet systems
over dyadic simplicial partitions of $I^d$ are a candidate for any dimension $d>1$.
For $m=3$ and $d=2$, the $12$-split Powell-Sabin spline spaces over dyadic triangulations of $I^2$ may potentially lead to such a construction, see \cite{Os1994}. However, we doubt that covering the range (\ref{PR6}) has any merit beyond academic interest.

\subsection{The exceptional case $s=0$, $1\le p < \infty$}\label{sec53} 
With a few exceptions, the literature about Besov function spaces such as $B_{p,q,1}^s$ deals only with parameters $s>0$ although the definition
given in Section \ref{sec22} formally leads to meaningful spaces in the case $s=0$ as well.
As our proof of Theorem \ref{theo2} reveals, the statement
$$
(B_{p,q,1}^s)'=\{0\}
$$
about the triviality of the dual space remains true also for $s=0$ if we are in the range $0<p<1$, $0<q<\infty$.
Thus, the question about the Schauder basis property of $H^d$ in $B_{p,q,1}^0$ is of interest 
only if $1\le p<\infty$. Since $H^d$ is a Schauder basis in $L_p=B_{p,\infty,1}^s$ for this range
(even unconditional if $1<p<\infty$), one would expect positive results also for $B_{p,q,1}^0$ in the case $0<q<\infty$. Unfortunately, our proof of the unconditionality results formulated in Theorem \ref{theo1} and sketched at the end of Section \ref{sec32} uses the assumption $s>0$ in an essential way.
We thus use an alternative argument inspired by \cite{Kr1978} to establish the following
result.
\begin{theo}\label{theo5} The Haar wavelet system $H^d$ is an unconditional Schauder basis
in $B_{p,q,1}^0(I^d)$ if and only if $1<p<\infty$, $0<q<\infty$. It is still a conditional
Schauder basis if $p=1$, $0<q<\infty$ (assuming blockwise ordering).
\end{theo}
{\bf Proof}. The proof uses the following well-known facts. Since $H^d$ is a Schauder basis
in $L_p$, $1\le p<\infty$, we have 
\be\label{ES}
E_k(g)_p \approx \|g-P_kg\|_{L_p}, \qquad k=0,1,\ldots,
\ee
for all $g\in L_p$, $1\le p<\infty$. Indeed, 
$$
E_k(g)_p \le \|g-P_kg\|_{L_p} \le \|g-s_k\|_{L_p}+\|P_k(g-s_k)\|_{L_p} \le CE_k(g)_p
$$
since $P_ks_k=s_k$ for the best approximating element $s_k\in S_k^d$ and since the 
partial sum operators $P_k$ are uniformly bounded in $L_p$, see Lemma \ref{lem3}.
Thus, if $P$ is a partial sum operator for $H^d$ of the form (\ref{defPS}) then
$E_l(Pg)_p=0$ for $l>k$ since $Pg\in S_{k+1}^d$, and
\be\label{ES1}
E_l(Pg)\le \|Pg-P_lPg\|_{L_p} = \|P(g-P_l)g\|_{L_p}\le C\|g-P_lg\|_{L_p}\le C E_l(g)_p
\ee
for $l=0,\ldots,k$, since $P$ and $P_l$ commute, and the $P$ are uniformly bounded in $L_p$. This, together with Lemma \ref{lem1} for $s=0$, establishes the uniform 
boundedness of the operators $P$, and consequently the Schauder basis property of $H^d$ in the
assumed natural ordering, in $B_{p,q,1}^0$ for all $1\le p<\infty$, $0<q<\infty$.

To prove the stronger unconditionality statement for $1<p<\infty$, we again resort to the
criterion for unconditionality stated in Lemma \ref{lem3}. Since $H^d$ is an unconditional
Schauder basis in  $L_p$, $1<p<\infty$, we have
\be\label{UncLp}
\|P_J g\|_{L_p}\le C \|g\|_{L_p},\qquad P_J g:=\sum_{h\in J} \theta_h\lambda_h(g) h,
\ee
for all $g\in L_p$ and all finite subsets $J$ of $H^d$. Since the projectors $P_J$ and $P_k$ commute,
as in (\ref{ES1}) this also implies 
\be\label{ES2}
E_k(P_J g)\le C E_k(g),\qquad k=0,1,\ldots, \qquad g\in L_p.
\ee 
We use (\ref{UncLp}-\ref{ES2}) in conjunction with Lemma \ref{lem1}, and bound the $B_{p,q,1}^0$ quasi-norm of $P_J g$ as follows:
\bea
\|P_J g\|_{B_{p,q,1}^0}^q &\le & C \|P_J g\|_{A_{p,q,1}^0}^q
= C\left(\|P_J g\|_{L_p}^q + \sum_{k=0}^\infty E_k(P_J g)_p^q\right)\\
&\le&  C\left(\|g\|_{L_p}^q + \sum_{k=0}^\infty E_k(g)_p^q\right)
= C \| g\|_{A_{p,q,1}^0}^q \le C\| g\|_{B_{p,q,1}^0}^q.
\eea
This proves the unconditionality of $H^d$ in $B_{p,q,1}^0$ for $1<p<\infty$, $0<q<\infty$.

Finally, to show that $H^d$ is not unconditional in $B_{p,q,1}^0$ if $p=1$, $0<q<\infty$,
we reuse the example $f_m$ from Section \ref{sec42}. Since Lemma \ref{lem2} also holds for $p=1$,
we compute
$$
E_l(f_m)_1 =\left\{ \ba{ll} 1,& l<m,\\ 0,& l\ge m,\ea  \right.
$$
and 
$$
E_l(g_{2k})_1 \left\{ \ba{ll} \approx 2k-l,& l<2k,\\  =0,& l\ge 2k,\ea  \right.
$$
where $g_{2k}=P_{J_k}f_m$ for a specific choice of $J_k\subset H^d$ and $2k\le m$.
Substitution into the formulas for the equivalent $A_{1,q,1}^0$ quasi-norms of these functions gives
$$
\|f_m\|_{B_{1,q,1}^0}^q\approx m,\qquad \|g_{2k}\|_{B_{1,q,1}^0}^q=\|P_{J_k}f_m\|_{B_{1,q,1}^0}^q
\approx \sum_{l<2k} (2k-l)^q \approx k^{q+1}.
$$
Now choose $k\approx m/2$ and let $m\to \infty$, to show that the operators $P_J$, $J\subset H^d$, are not uniformly
bounded in $B_{1,q,1}^0$, $0<q<\infty$. Consequently, $H^d$ is only a conditional Schauder basis in
$B_{1,q,1}^0$.
\hfill $\Box$

\medskip
One may ask if an analog of Theorem \ref{theo3}, i.e., an isomorphism with a certain sequence space
in terms of Haar coefficients, holds for $B_{p,q,1}^0$, $1<p<\infty$, $0<q<\infty$, as well. 
In principle, the answer is yes but it looks more complicated than the description via
weighted $\ell_q(\ell_p)$ spaces. Indeed, the Littlewood-Paley-type characterization for $L_p$ in terms of $H^d$ expansions,
namely the norm equivalence
\be\label{LiPa}
\|f\|_{L_p} \approx \| (\sum_{l=0}^\infty \sum_{h\in H^d_l} 2^{-ld}\lambda_h(f)^2 |h|^2)^{1/2}\|_{L_p},\quad f\in L_p,\quad 1<p<\infty,
\ee
see \cite[Corollary 2.28, (2.223)]{Tr2010}, can be used in conjunction with (\ref{ES}) and Lemma \ref{lem1}:
\bea
\|g\|_{B_{p,q,1}^0}^q &\approx& \|g\|_{A_{p,q,1}^0}^q
\approx \|g\|_{L_p}^q + \sum_{k=0}^\infty \|g-P_kg\|_{L_p}^q\\
&\approx& \sum_{k=0}^\infty \left\|\left(\sum_{l=k}^\infty\sum_{h\in H^d_l} 2^{-ld}\lambda_h(f)^2 |h|^2\right)^{1/2}\right\|_{L_p}^q.
\eea
Recall that $|h|^2=\chi_\Delta$, where $\Delta$ is the support cube of the Haar function $h\in H^d$.

We doubt that the above equivalent quasi-norm for $B_{p,q,1}^0$, $1<p<\infty$, $0<q<\infty$, can be made more explicit, except for
the special case $p=2$, where the orthogonality of $H^d$ leads to simplifications:
$$
\|g\|_{B_{2,q,1}^0}^q\approx
\sum_{k=0}^\infty \left(\sum_{l=k}^\infty \sum_{h\in H^d_l} 2^{-ld}\lambda_h(g)^2\right)^{q/2}, \quad 0<q<\infty.
$$
This further simplifies if we also set $q=2$:
$$
\|g\|_{B_{2,2,1}^0}^2\approx
\sum_{k=0}^\infty \sum_{h\in H^d_k} (k+1)2^{-kd}\lambda_h(g)^2, \quad 0<q<\infty.
$$

\subsection{Tensor-product Haar system $\tH^d$ as basis in $B_{p,q,1}^s$}\label{sec52} 
The historically first constructions of Schauder bases for Banach spaces
of smooth functions over higher-dimensional cubes and manifolds \cite{CiDo1972,CiFi1982,CiFi1983a} exclusively used tensor products of univariate Schauder bases. Later, due to the desire to work with systems with better support localization and the need to cover quasi-Banach spaces as well, wavelet-type constructions became more popular. As a matter of fact, in more recent texts such as \cite{Tr2010}, the tensor-product construction is examined only with respect to function spaces of dominating mixed smoothness. 

This raises the following question: Can the Haar tensor-product system $\tH^d$ be a Schauder basis in any of the separable Besov function spaces $B_{p,q,1}^s$ it is contained in? Note that if we write
$$
\tH^d = \cup_{k=0}^\infty \tH_k^d, 
$$
where $\tH^d_k=\{\th\in S^d_k:\;\th\not\in S_{k-1}^d$ for $k\ge 1$ and $\tH_0^d=H^d_0$,
and order by blocks, then
the set of partial sum operators $\tP$ associated with $\tH^d$ expansions contains 
$\{P_k\}_{k\in \mathbb{Z}_+}$ as a subsequence. Consequently, such an ordering of $\tH^d$ assumed, $\{\tP\}$ is uniformly bounded in $B_{p,q,1}^s$ only if the Haar wavelet system $H^d$ with the natural ordering is a Schauder basis in this space. When this is the case, one may hope that $\tH^d$ is a Schauder basis as well.
	
For $1\le p <\infty$ this is indeed the case: The properly ordered tensor-product Haar system
$\tH^d$ is a Schauder basis in $B_{p,q,s}^s$, $1\le p<\infty$, $0<q<\infty$, $0\le s <1/p$, and is even unconditional if $1<p<\infty$. We believe that this is known but could not locate the corresponding statement in the literature. 	
However, for $0<p<1$ a simple example shows that $\tH^d$ is not a Schauder basis in $B_{p,q,s}^s$, independently of the ordering in $\tH^d$.
\begin{theo}\label{theo4} Let $d>1$, $0<q<\infty$, $0\le s <1/p$.\\
a) The tensor-product Haar system $\tH^d$ is an unconditional Schauder basis in $B_{p,q,1}^s(I^d)$ if and only if $1<p<\infty$,
for $p=1$ it remains a conditional Schauder basis for a particular blockwise ordering.\\	
b) If $0<p<1$ then the family of $L_2$ orthogonal projectors 
	$$
	\Lambda_{\th}:\, f \to \Lambda_{\th} f := \|\th\|_{L_2}^{-2} \left(\int_{I^d} f\th \,dx\right) \th, \qquad \th\in \tH^d,
	$$
	onto one-dimensional subspaces spanned by tensor-product Haar functions $\th\in \tH^d$
	cannot be uniformly bounded in $B_{p,q,1}^s(I^d)$.
	As a consequence, $\tH^d$ cannot be a Schauder basis in $B_{p,q,1}^s(I^d)$ for $0<p<1$,
	independently of the ordering within $\tH^d$.
\end{theo}

{\bf Proof}. 
\emph{Step1}. We start with part b). For $k=1,2,\ldots$, consider the characteristic function $f_k:=\chi_{[0,2^{-k})^d}\in S^d_k$ and
the particular tensor-product Haar function
$$
\th_k:= h_{[0,2^{-(k-1)})}\otimes \chi_I\otimes \ldots \otimes \chi_I =
h_{2^{k-1}+1}  \otimes h_1\otimes \ldots \otimes h_1\in \tH^d_k,
$$
see the notation in Section \ref{sec21}. Using Lemma \ref{lem1} and \ref{lem2}, we have
$$
\|f_k\|_{B_{p,q,1}^s}^q\approx \|f_k\|_{L_p}^q +\sum_{l=0}^{k-1} 2^{lsq}E_l(f_k)_p^q =
\|f_k\|_{L_p}^q(1+\sum_{l=0}^{k-1} 2^{lsq})\approx 2^{-kdq/p}a_{q,s,k},
$$
where $a_{q,s,k}:=\sum_{l=0}^{k} 2^{lsq}$, and $0\le s<1/p$, $0<q<\infty$.
Similarly,
$$
\|\th_k\|_{B_{p,q,1}^s}^q\approx \|\th_k\|_{L_p}^q +\sum_{l=0}^{k-1} 2^{lsq}E_l(\th_k)_p^q \approx
\|\th_k\|_{L_p}^q(1+\sum_{l=0}^{k-1} 2^{lsq})\approx 2^{kq/p}a_{q,s,k},
$$
where, by applying Lemma \ref{lem2} locally on cubes $\Delta\in T_l^{d}$
with $l<k$, we could use that
$E_l(\th_k)_p = \|\th_k\|_{L_p} = 2^{-(k-1)/p)}$ for all $l<k-1$ and
$E_{k-1}(\th_k)_p= 2\cdot 2^{-k/p} = 2^{1-1/p}\|\th_k\|_{L_p}$ for $l=k-1$.

Putting things together, with
$$
\|\th_k\|_{L_2}^2 = 2^{-k+1}, \qquad \int_{I^d} f_k\th_k \,dx = 2^{-kd},
$$
and
$$
\|\Lambda_{\th_k}f_k\|_{B_{p,q,1}^s}^q = (2^{k-1} \cdot 2^{-kd})^q \|\th_k\|_{B_{p,q,1}^s}^q,
$$
we arrive at
$$
\|\Lambda_{\th_k}\|_{B_{p,q,1}^s\to B_{p,q,1}^s}^q\ge \frac{\|\Lambda_{\th_k}f_k\|_{B_{p,q,1}^s}^q}{\|f_k\|_{B_{p,q,1}^s}^q}\approx 
\frac{2^{-k(d-1)q}2^{-kq/p}}{2^{-kdq/p}}=2^{k(1/p-1)(d-1)q}.
$$
Since $0<p<1$ and $d>1$ this shows that 
$\{\Lambda_{\th_k}\}_{k=1,2,\ldots}$ and, consequently, the whole family $\{\Lambda_{\th}\}_{\th\in\tH^d}$ cannot be uniformly bounded on $B_{p,q,1}^s$.

Finally, since $\Lambda_{\th}$ is the difference of two consecutive partial sum operators
for the series expansion with respect to $\tH^d$, independently of the ordering within $\tH^d$,
the sequence of partial sum operators cannot be uniformly bounded on $B_{p,q,1}^s$ either.
Due to Lemma \ref{lem3} this shows that $\tH^d$ cannot have the Schauder basis property
in $B_{p,q,1}^s$ for the indicated parameter range if $0<p<1$.
\hfill $\Box$

\smallskip
\emph{Step 2}. We next prove the unconditionality of $\tH^d$ 
in $B_{p,q,1}^s$ for all parameters       
\be\label{PR7}
1 \le  p <\infty, \quad 0<q<\infty, \quad 0\le s <1/p.
\ee
To this end, recall from Lemma \ref{lem1} and (\ref{ES}) that
\be\label{ABP}
\|f\|_{B_{p,q,1}^s}^q\approx \|f\|_{A_{p,q,1}^s}^q \approx \|f\|_{L_p}^q + \sum_{k=0}^\infty
(2^{ks}\|f-P_kf\|_{L_p})^q, \qquad f\in B_{p,q,1}^s,
\ee
for the parameters in (\ref{PR7}). Using  
$$
\|f-P_k\|_{L_p} \le \sum_{l=k}^\infty \|(P_{l+1}-P_l)f\|_{L_p}
$$
and (\ref{Hardy}), this further transforms to
\be\label{ABP1}
\|f\|_{B_{p,q,1}^s}^q\approx \sum_{k=0}^\infty
(2^{ks}\|P_k-P_{k-1}f\|_{L_p})^q, \qquad f\in B_{p,q,1}^s,
\ee
if $ 0<s<1/p$ (we have set $P_{-1}=0$).

Similar to (\ref{LiPa}) for $H^d$, we have a Littlewood-Paley-type norm equivalence 
\be\label{NELp}
\|g\|_{L_p} \approx \left\|\left(\sum_{\th\in\tH^d} (\mu(\mathrm{supp}(\th))^{-1}\lambda_{\th}(g)\th)^2\right)^{1/2}\right\|_{L_p}, \qquad g\in L_p,\quad 1<p<\infty,
\ee
for $\tH^d$, where $\mu(\cdot)$ denotes the Lebesgue measure in $I^d$, see \cite[Corollary 2.28, (2.223)]{Tr2010}.
Substituted into (\ref{ABP1}), this yields the following characterization of the 
$B_{p,q,1}^s$ quasi-norm if $1<p<\infty$, $0<q<\infty$, and $0<s<1/p$:
$$%\be\label{ABP2}
\|f\|_{B_{p,q,1}^s}^q\approx \sum_{k=0}^\infty
\left(2^{ks}\left \| \left(\sum_{\th\in\tH^d_k} (\mu(\mathrm{supp}(\th))^{-1}\lambda_{\th}(g)\th)^2\right)^{1/2} \right\|_{L_p}\right)^q, \qquad f\in B_{p,q,1}^s.
$$%\ee
For $s=0$, we obtain a similar characterization of the 
$B_{p,q,1}^0$ quasi-norm if we substitute (\ref{NELp}) into (\ref{ABP}).
This shows that $\tH^d$ is an unconditional Schauder basis in $B_{p,q,1}^s$ for the parameters in (\ref{PR7}) if $1<p<\infty$. 

\smallskip
\emph{Step 3}.
For $p=1$ we only establish the Schauder basis property for a specific blockwise ordering of $\tH^d$. The idea is old, see e.g. \cite[Section 11]{CiDo1972} or \cite[Section 4-5]{CiFi1982}, in the latter paper, the Haar case corresponds to the parameter settings $r=1$, $k=0$. For clarity, we restrict ourselves to $d=2$, the general case can be handled by induction in $d$. Consider the natural ordering $H=\{h_n\}_{n\in \mathbb{N}}$ of the univariate Haar system $H$ introduced in Section \ref{sec21}, and denote the partial sum operators by $Q_n$. These operators are $L_2$ ortho-projectors onto $\mathrm{span}(\{h_l\}_{l=1,\ldots,n})$, and uniformly bounded on $L_1$.  
	
With this ordering for $H$, we enumerate the blocks $\tH_{k}^2$ of the tensor-product Haar system
$$
\tH^2=\{\th_{i_1,i_2}:=h_{i_1}\otimes h_{i_2},\quad (i_1,i_2)\in\mathbb{N}^2\}
$$
as follows: We set $\tH^2_0=\{\th_1:=\th_{1,1}\}$ for $k=0$, and 
\be\label{DHk2}
\tH^2_{k+1}=\underbrace{\left(\cup_{i=1}^{2^{k}}\{\th_{2^{k}+i,n}\}_{n=1,\ldots,2^{k+1}}\right)}_{\tH'_{k+1}} \cup 
\underbrace{\left(\cup_{i=1}^{2^{k}}\{\th_{n,2^{k}+i}\}_{n=1,\ldots,2^{k}}\right)}_{\tH''_{k+1}}
\ee
for $k\ge 0$. The enumeration of the functions within the subblocks $\tH'_{k+1}$ and $\tH''_{k+1}$ in (\ref{DHk2}) is assumed lexicographic with respect to the index pairs $(i,n)$. With this at hand, any partial sum operator $\tP$  with respect to $\tH^2$ is the linear combination of a few projectors $Q_{n_1,n_2}=Q_{n_1}\otimes Q_{n_2}$.  Indeed, similar to
(\ref{defPS}), any partial sum operator $\tP$ takes the form
$$%\be\label{defPS1}
\tP g=P_kg + \Delta\tP g;\qquad \Delta \tP g:=\sum_{\th\in \bar{\tH}^2_{k+1}} \lambda_{\th}(g)\th \in S_{k+1}^2.
$$%\ee
for some $k=0,1,\ldots$ and some section $\bar{\tH}^2_{k+1}$ of ${\tH}^2_{k+1}$ taken in the described order. Obviously, $P_k=Q_{2^k,2^k}$. If the section $\bar{\tH}^2_{k+1}$ is contained in $\tH'_{k+1}$ then
\bea
\Delta \tP &=& (Q_{2^k+i-1}-Q_{2^k})\otimes Q_{2^{k+1}} + (Q_{2^k+i}-Q_{2^k+i-1})\otimes Q_n\\
&=& Q_{2^k+i-1,2^{k+1}} - Q_{2^k,2^{k+1}} + Q_{2^k+i,n} - Q_{2^k+i-1,n}
\eea
for some $n=1,\ldots, 2^{k+1}$ and  $i=1,\ldots, 2^k$. Otherwise, we have
\bea
\Delta \tP &=& (Q_{2^{k+1}}-Q_{2^k})\otimes Q_{2^{k+1}} + Q_{2^k}\otimes (Q_{2^k+i-1}-Q_{2^k}) + Q_n\otimes(Q_{2^k+i}-Q_{2^k+i-1})\\
&=& Q_{2^{k+1},2^{k+1}}-Q_{2^{k},2^{k+1}} +Q_{2^{k+},2^{k}+i-1}-Q_{2^{k},2^{k}}+
Q_{n,2^{k}+i}-Q_{n,2^{k}+i-1},
\eea
for some $n=1,\ldots, 2^{k}$ and  $i=1,\ldots, 2^k$. Altogether, this shows that $\tP$ is always a linear combination of a few tensor-product projectors $Q_{n_1,n_2}$. Since tensor-products of uniformly $L_1$ bounded operators are uniformly $L_1$ bounded as well, it follows that
\be\label{L1bound}
\|\tP g\|_{L_1} \le C\|g\|,\qquad g\in L_1,
\ee
uniformly for all partial sum operators $\tP$. Thus, $\tH^d$ is a Schauder basis in $L_1$ for $d=2$. By induction, this holds for all $d>1$. This is most probably known, and has been proved here only in the absence of a proper reference.

For the parameters under consideration, the Schauder basis property of $\tH^d$ in $B_{1,q,1}^s$ follows now from (\ref{L1bound}) using the same arguments as in the proof of Theorem \ref{theo5} in Section \ref{sec53}. 

\smallskip
\emph{Step 4}.
Finally, it is easy to see that $\tH^d$ is not unconditional in $B_{1,q,1}^s$ for 
any $0\le s <1$, $0<q<\infty$. If it were, then any finite subset of $\tH^d$ must be uniformly unconditional as well. But this is not the case. Consider the finite subsets
$$
J_k:=\{h_n\otimes h_{2^{k}+1}\otimes \ldots\otimes  h_{2^k+1}\}_{n=1,\ldots,2^{k+1}}\subset \tH^d_{k+1},\qquad k=0,1,\ldots.
$$
This $J_k$ is a finite section of $H\otimes h_{2^k+1}\otimes \ldots\otimes  h_{2^k+1}$. Since $H$ is not unconditional in $L_1(I)$, the subsets $J_k$ cannot be uniformly unconditional in $L_1$ either. This means that there is a sequence of subsets $J'_k\subset J_k\subset \tH_{k+1}^d$, and
a sequence of functions
$$
g_k\in \mathrm{span}(\tH^d_{k+1})=\mathrm{span}(H^d_{k+1})=S_{k+1}^d\ominus_{L_2} S_k^d,\qquad k=0,1,\ldots,
$$
such that
\be\label{UncondTP}
\frac{\|\tP_{J'_k}g_k\|_{L_1}}{\|g_k\|_{L_1}} \to \infty, \qquad k\to \infty,\qquad 
\ee 
where $\tP_{J'_k}$ denotes the partial sum projector with respect to the subset
$J'_k\subset \tH^d_{k+1}$.

Consider now the $B_{1,q,1}^s$ quasi-norms of $g_k$ and $\tP_{J'_k}g_k$. Since these functions belong to $S^d_{k+1}$ but are $L_2$-orthogonal to $S_k^d$ we have 
$E_l(g_k)_1=E_l(\tP_{J'_k}g_k)_1=0$ for $l>k$ while for $l\le k$ we have $P_lg_k=P_l\tP_{J'_k}g_k=0$ and according to (\ref{ES}) 
$$
E_l(g_k)_{1}\approx \|g_k\|_1,\qquad E_l(\tP_{J'_k}g_k)_1\approx \|\tP_{J'_k}g_k\|_1.
$$
If we substitute this into the expressions for the $A_{1,q,1}^s$ norms then using Lemma \ref{lem1} we get 
$$
\|\tP_{J'_k}\|_{B_{1,q,1}^s\to B_{1,q,1}^s}
\ge\frac{\|\tP_{J'_k}g_k\|_{B_{1,q,1}^s}}{\|g_k\|_{B_{1,q,1}^s}} 
\approx \frac{\|\tP_{J'_k}g_k\|_{A_{1,q,1}^s}}{\|g_k\|_{A_{1,q,1}^s}} 
\approx \frac{\|\tP_{J'_k}g_k\|_{L_1}}{\|g_k\|_{L_1}}.
$$
Due to (\ref{UncondTP}), this contradicts the unconditionality criterion formulated in Lemma \ref{lem3}. \hfill $\Box$

\subsection{The exceptional case $q=\infty$}\label{sec54} 
Since $B_{p,\infty,1}^s$ is non-separable for $0<s\le \max(1,1/p)$, it cannot possess Schauder bases. %since the span of a countable system of functions cannot be dense in a non-separable quasi-Banach space. 
However, one can talk about
the (unconditional) Schauder basis property of a system $(f_m)_{m\in\mathbb{N}}$ in a quasi-Banach space $X$  with respect to its closure in $X$. In this case, $(f_m)_{m\in\mathbb{N}}$ is called (unconditional) basis sequence in $X$. In particular, one can ask if $H^d$ is an (unconditional) basis sequence
in  $B_{p,\infty,1}^s(I^d)$ if $0<s\le 1/p$.
% (or, equivalently, if $H^d$ is an (unconditional) Schauder basis in the space
%$$
%b_{p,\infty,1}^s = \{f\in L_p:\, \omega(t,f)_p=\mathrm{o}(t^s),\;t\to 0\}, \quad 0<s\le 1/p,
%$$
%which, equipped with the $B_{p,\infty,1}^s$ quasi-norm, is the closure of $H^d$ in $B_{p,\infty,1}^s$). 
Since $B_{p,\infty,1}^0=L_p$, the case $s=0$ does not need consideration.

For $d=1$ and $1\le p <\infty$, a detailed study of the Schauder basis property and unconditionality
of the univariate  Haar system $H$ has been provided by Krotov \cite{Kr1978,Kr1982}
for the generalized Nikolski classes 
$$
\Lambda^\omega_p(I):=\{f\in L_p(I):\, \|f\|_{\Lambda^\omega_p}:=\|f\|_{L_p} +\sup_{0<t<1}
\frac{\omega(t,f)_p}{\omega(t)} < \infty\},
$$
where $\omega(t)$ is an appropriate comparison function. If $\omega(t)=t^s$, we have $\Lambda^\omega_p(I)=
B_{p,\infty,1}^s(I)$ as partial case. In particular, the results in \cite{Kr1978,Kr1982} for $1\le p<\infty$
imply that $H$ is an unconditional basis sequence in $B_{p,\infty,1}^s(I)$ for $0<s< 1/p$, while for $s=1/p$  it is a conditional basis sequence.  More recently, \cite{GSU2019a} studied basis sequence
properties of Haar wavelet systems in distributional Besov spaces $B_{p,\infty}^s(\mathbb{R}^d)$ and
$B_{p,\infty}^s$ for $0<p<\infty$. 

%To prove the basis sequence property for $H^d$, uniform bounds for the partial sum operators $P$ and $P_J$ in $B_{p,\infty,1}^s$ are  necessary and sufficient. 
Implicitly, our considerations in Section \ref{sec3} and \ref{sec4}
already answer the questions about the basis sequence properties for $H^d$ in $B_{p,\infty,1}^s$. Without proof, we state the results. The examples in Section \ref{sec41} show that for $0<s\le d(1/p-1)$,
$0<p<1$, bounded linear functionals on $B_{p,\infty,1}^s$ must be trivial on the closure of $H^d$. Consequently, $H^d$ is then not a basis sequence in $B_{p,\infty,1}^s$. For the parameter range
$$
\max(d(1/p-1),0) < s < 1/p, \qquad (d-1)/d < p < \infty,
$$
in agreement with the above cited results for $d=1$ and by slightly modifying the arguments for Theorem \ref{theo3}, we can prove that $H^d$ is an unconditional basis sequence in $B_{p,\infty,1}^s$. More precisely, the map $\Lambda$ defined in (\ref{LambdaMap}) provides an isomorphism between $B_{p,\infty,1}^s$ and a closed subspace of the weighted $\ell_\infty(\ell_p))$ sequence space consisting of all sequences $\Gamma=(\gamma_h)_{h\in H^d}$ with finite quasi-norm
$$%\be\label{Infty}
\|\Gamma\|_{\ell_\infty(\ell_p)} :=\sup_{k\ge 0} 2^{k(s-d/p)}\left(\sum_{h\in H^d_k} |\gamma_h|^p \right)^{1/p}.
$$%\ee
This subspace is characterized by the additional condition
$$
2^{k(s-d/p)}\left(\sum_{h\in H^d_k} |\gamma_h|^p \right)^{1/p} \to 0,\qquad k\to\infty.
$$
With the exception of the parameter range $1\le s <1/p$ and $(d-1)/d <p<1$, this result can also be recovered from the literature, see, e.g., \cite{Ci1975,Ro1976,Tr2010}. 

For $s=1/p$, $(d-1)/d < p < \infty$, the Haar wavelet system $H^d$ is only a conditional basis sequence. Since $B^{1/p}_{p,\infty,1}$ cannot be characterized in terms of the sequence $(E_k(f)_p)_{k\in\mathbb{Z}_+}$, and one needs to work with the original definition
of the $B^{1/p}_{p,\infty,1}$ quasi-norm in terms of moduli of smoothness. 
We leave it to the reader to fill in the details.

%\section*{References} 

\end{document}